\input amstex
\documentstyle{amsppt}

\input label.def
\input degt.def

{\catcode`\@11
\gdef\proclaimfont@{\sl}}

\newtheorem{Remark}\Remark\thm\endAmSdef
\conjecture\thm\endproclaim

\input epsf
\def\picture#1{\epsffile{#1-bb.eps}}
\def\plot#1{\centerline{\epsfxsize.8\hsize\picture{plot-#1}}}

\def\dash{\item"\hfill--\hfill"}
\def\Dashes{\widestnumber\item{--}\roster}
\def\endDashes{\endroster}

\def\ie{\emph{i.e.}}
\def\eg{\emph{e.g.}}
\def\cf{\emph{cf}}

\def\etc{\emph{etc}}

\let\Ga\alpha
\let\Gb\beta
\let\Gg\gamma
\let\Gd\delta
\let\Gs\sigma
\let\Gf\varphi
\let\Gr\rho
\def\bGa{\bar\Ga}
\def\bGb{\bar\Gb}
\def\bGg{\bar\Gg}

\def\1{^{-1}}
\def\2{^{-2}}

\def\xt{x_t}
\def\yt{y_t}

\def\bpi{\bar\pi_1}

\loadbold
\def\bA{\bold A}
\def\bD{\bold D}
\def\bE{\bold E}

\let\splus\oplus
\let\bigsplus\bigoplus
\let\onto\twoheadrightarrow

\def\Conv{\operatorname{conv}}
\def\Aut{\operatorname{Aut}}

\def\CG#1{\Z_{#1}}     
\def\BG#1{\Bbb B_{#1}} 
\def\SG#1{\Bbb S_{#1}} 

\def\Cp#1{\Bbb P^{#1}}
\def\Rp#1{\Bbb P^{#1}_\R}
\def\B{\bar B}
\def\L{\bar L}
\def\CM{\Cal M}
\def\tCM{\tilde\CM}

\def\MB{D} 
\def\disk{\Delta}

\def\PSL{\operatorname{\text{\sl PSL}}}
\def\SL{\operatorname{\text{\sl SL}}}
\def\gcd{\operatorname{g.c.d.}}
\def\<#1>{\langle#1\rangle}
\def\ls|#1|{\mathopen|#1\mathclose|}

{\let\1\gdef\let\+\relax
\gdef\setcat{\catcode`\a\active \catcode`\d\active \catcode`\e\active
\catcode`\+\active}
\setcat
\1a[#1]{\bA_{#1}} \1d[#1]{\bD_{#1}} \1e[#1]{\bE_{#1}}
}
{\catcode`\+\active\gdef+{\splus}}
\def\getline#1,#2[#3],#4[#5],#6,#7[#8]{#2&#1\cr}
{\obeylines%
\gdef\TAB#1{\vtop\bgroup\openup1pt\obeylines\let
\getline\setcat\expandafter\halign\expandafter\bgroup#1\cr}}
\def\ENDTAB{\crcr\egroup\egroup}

\def\NonTorus{$##$\hss&$##$\hss}

\def\Tab{\vtop\bgroup\openup1pt\halign\bgroup$##$\hss\cr}

\topmatter

\author
Alex Degtyarev
\endauthor

\title
Fundamental groups of symmetric sextics
\endtitle

\address
Department of Mathematics,
Bilkent University,
06800 Ankara, Turkey
\endaddress

\email
degt\@fen.bilkent.edu.tr
\endemail

\abstract
We study the moduli spaces and compute the fundamental groups of
plane sextics of torus type with at least two type~$\bold{E}_6$
singular points. As a simple
application, we compute the fundamental
groups of $125$ other sextics, most of which are new.
\endabstract



\endtopmatter

\document

\section{Introduction}

\subsection{Principal results}
Recall that a plane sextic~$B$ is said to be of \emph{torus type}
if its equation can be represented in the form $p^3+q^2=0$, where
$p$ and~$q$ are certain homogeneous polynomials of degree~$2$
and~$3$, respectively. Alternatively, $B\subset\Cp2$ is of torus
type if and only if it is the ramification locus of a projection
to~$\Cp2$ of a cubic surface in~$\Cp3$. A representation of the
equation in the form $p^3+q^2=0$ (up to the obvious equivalence)
is called a \emph{torus structure} of~$B$. A singular point~$P$
of~$B$ is called \emph{inner} (\emph{outer}) with respect to a
torus structure~$(p,q)$ if $P$ does (respectively, does not)
belong to the intersection of the conic $\{p=0\}$ and the cubic
$\{q=0\}$. The sextic~$B$ is called \emph{tame} if all its
singular points are inner. Note that, according
to~\cite{degt.Oka}, each sextic~$B$ considered in this paper has a
unique torus structure; hence, we can speak about inner and outer
singular points of~$B$.
For the reader's convenience, when
listing the set of singularities of a sextic of torus type,
we indicate the inner singularities by
enclosing them in parentheses.

Apparently, it was O.~Zariski~\cite{Zariski} who first understood
the importance of sextics of torus type. Since then, they have
been a subject of intensive study. For details and further
information,
we refer to
M.~Oka, D.~T.~Pho~\cite{OkaPho.moduli},
\cite{OkaPho} (topology, sets of singularities, moduli,
fundamental groups), H.~Tokunaga~\cite{Tokunaga} (algebro-geometric
approach), and A.~Degtyarev~\cite{degt.Oka}.

In recent paper~\cite{degt.8a2}, we described the moduli spaces
and calculated the fundamental groups of all sextics of torus type
of weight~$8$ and~$9$ (in a sense, those with the largest
fundamental groups). The approach used in~\cite{degt.8a2},
reducing sextics to maximal
trigonal curves, was also helpful in the study
of some other sextics with nonabelian groups
(see~\cite{degt.Oka3}), and then, in~\cite{symmetric}, we
classified all irreducible sextics for which this approach should
work. The purpose of this paper is to treat one of the classes
that appeared in~\cite{symmetric}: sextics with at least two
type~$\bE_6$ singular points; they are reduced to trigonal curves
with the set of singularities $\bE_6\splus\bA_2$.
Our principal results are
Theorems~\ref{th.moduli} and~\ref{th.group} below.

\midinsert
\table
Sextic with two type~$\bE_6$ singular points
\endtable\label{tab.list}
\hbox to\hsize{\hss\def\1{\llap{$^*$\,}}
\Tab
\1 (3\bE_6)\splus\bA_1\cr
   (3\bE_6)\cr
\1 (2\bE_6\splus\bA_5)\splus\bA_2\cr
   (2\bE_6\splus\bA_5)\splus\bA_1\cr
   (2\bE_6\splus\bA_5)
\ENDTAB\hss\Tab
\1 (2\bE_6\splus2\bA_2)\splus\bA_3\cr
\1 (2\bE_6\splus2\bA_2)\splus\bA_2\cr
\1 (2\bE_6\splus2\bA_2)\splus2\bA_1\cr
   (2\bE_6\splus2\bA_2)\splus\bA_1\cr
   (2\bE_6\splus2\bA_2)
\ENDTAB
\hss}
\endinsert

\theorem\label{th.moduli}
Any sextic of torus type with at least two
type~$\bE_6$ singular points has one of the sets of singularities
listed in Table~\ref{tab.list}. With the exception of
$(2\bE_6\splus\bA_5)\splus\bA_2$, the moduli space of sextics of
torus type realizing each set of singularities in the table is
rational \rom(in particular, it is nonempty and
connected\rom)\rom; the moduli space of sextics with the set of
singularities $(2\bE_6\splus\bA_5)\splus\bA_2$ consists of two
isolated points, both of torus type.
\endtheorem

Note
that we do not assume \emph{a priori} that the curves are
irreducible or have simple singularities only.
Both assertions hold automatically for any sextic with
at least two type~$\bE_6$ singular points, see the beginning of
Section~\ref{proof.moduli}.

Theorem~\ref{th.moduli} is proved in Section~\ref{proof.moduli}.
The two classes of sextics realizing the
set of singularities $(2\bE_6\splus\bA_5)\splus\bA_2$ were first
discovered in Oka, Pho~\cite{OkaPho.moduli}.
The sets of singularities that can be
realized by sextics of torus type are
also listed in~\cite{OkaPho.moduli}.
Note
that the list
given by Table~\ref{tab.list}
can also be obtained from the results of J.-G.~Yang~\cite{Yang},
using the characterization of irreducible
sextics of torus type found
in~\cite{degt.Oka}. The deformation classification can be obtained
using~\cite{JAG}.

\Remark
A simple calculation using~\cite{JAG} or~\cite{Yang}
and the characterization of irreducible
sextics of torus type found
in~\cite{degt.Oka}
shows that the sets of
singularities marked with a $^*$ in Table~\ref{tab.list} are realized by
sextics of torus type only. Each of the
remaining
five sets of
singularities is also realized by a single deformation family of
sextics not of torus type, see A.~\"Ozg\"uner~\cite{Aysegul} for
details. Furthermore, Table~\ref{tab.list} lists all sets of
singularities
of plane sextics, both of and not of torus type,
containing at least two
type~$\bE_6$ points.
\endRemark

\theorem\label{th.group}
Let $B$ be a sextic of torus type whose set of
singularities~$\Sigma$ is one of those listed in
Table~\ref{tab.list}. Then the fundamental group
$\pi_1:=\pi_1(\Cp2\sminus B)$ is as follows\rom:
\roster
\item\local{G3.group}
if $\Sigma=(2\bE_6\splus2\bA_2)\splus\bA_3$, then $\pi_1$ is the
group~$G_3$ given by~\eqref{eq.G3}\rom;
\item\local{G0.group}
if $\Sigma=(3\bE_6)\splus\bA_1$ or $(2\bE_6\splus2\bA_2)\splus2\bA_1$,
then $\pi_1=G_0:=\BG4/\Gs_2\Gs_1^2\Gs_2\Gs_3^2$\rom;
\item\local{two.groups}
if $\Sigma=(2\bE_6\splus\bA_5)\splus\bA_2$, then, depending on the
family, $\pi_1$ is one of
the groups~$G_2'$, $G_2''$ given by~\eqref{eq.G2.1}
and~\eqref{eq.G2.2}, respectively\rom;
\item\local{braid.group}
otherwise, $\pi_1=\BG3/(\Gs_1\Gs_2)^3$.
\endroster
\rom(Here, $\{\Gs_1,\ldots,\Gs_{n-1}\}$ is a canonical basis for
the braid group~$\BG{n}$ on $n$ strings.\rom)
\endtheorem

The fundamental groups are calculated in \S\ref{S.group}. An
alternative presentation of the groups~$G_2'$, $G_2''$
mentioned in~\iref{th.group}{two.groups} is found in C.~Eyral,
M.~Oka~\cite{EyralOka}, where it is conjectured that
the two groups are not isomorphic. We suggest to attack this
problem studying the relation between~$G_2'$, $G_2''$ and the local
fundamental group at the type~$\bA_5$ singular point, \cf.
Proposition~\ref{a5->>pi} and Conjecture~\ref{a5not->>pi}. The
group of a sextic of torus type with the set of singularities
$(2\bE_6\splus\bA_5)\splus\bA_1$, see~\iref{th.group}{braid.group},
is also found in~\cite{EyralOka}; the group of a sextic with the
set of singularities $(3\bE_6)\splus\bA_1$,
see~\iref{th.group}{G0.group}, as well as the groups of the three
tame sextics listed in Table~\ref{tab.list} (the sets of
singularities $(3\bE_6)$, $(2\bE_6\splus\bA_5)$, and
$(2\bE_6\splus2\bA_2)$)
are found in Oka, Pho~\cite{OkaPho}.

With the possible exception of $G_2'$, $G_2''$, all groups listed in
Theorem~\ref{th.group} are `geometrically' distinct in the
sense of the following theorem.

\theorem\label{th.proper}
All epimorphisms
$$
G_3\onto G_0\onto\BG3/(\Gs_1\Gs_2)^3,\qquad
G_2',G_2''\onto\BG3/(\Gs_1\Gs_2)^3
$$
induced by the respective
perturbations of the curves \rom(\cf. Zariski~\cite{Zariski}\rom)
are proper, \ie, they are not isomorphisms.
\endtheorem

This theorem is proved in Section~\ref{proof.proper}. Some of the
statements follow from the previous results by Eyral,
Oka~\cite{EyralOka} and Oka, Pho~\cite{OkaPho}.

As a further application of Theorem~\ref{th.group},
we use the presentations obtained and the
results of~\cite{degt.8a2} to compute the fundamental groups of
eight sextics of torus type and $117$ sextics not of torus type
that are not covered by M.~V.~Nori's theorem~\cite{Nori}, see
Theorems~\ref{th.nontorus} and~\ref{th.torus}. As for most
sets of singularities the connectedness of the moduli space has
not been established (although expected),
we state these results in the form of existence.

\subsection{Contents of the paper}
In~\S\ref{S.model}, we use the results of~\cite{symmetric} and
construct the \emph{trigonal models} of sextics in question, which
are pairs $(\B,\L)$, where $\B$ is a (fixed) trigonal curve
in the Hirzebruch surface~$\Sigma_2$ and $\L$ is a (variable)
section. We study the conditions on~$\L$ resulting in a particular
set of singularities of the sextic. As a consequence, we obtain
explicit equations of the sextics and rational parameterizations of
the moduli spaces. Theorem~\ref{th.moduli} is proved here.

In~\S\ref{S.vanKampen}, we present the classical Zariski--van
Kampen method~\cite{vanKampen} in a form suitable for curves on
Hirzebruch surfaces.
The contents of this section is a formal account of a few
observations found in~\cite{degt.Oka3}
and~\cite{degt.kplets}.

In~\S\ref{S.group}, we apply the classical Zariski--van Kampen
theorem to the trigonal models constructed above and obtain
presentations of the fundamental groups. The main advantage of
this approach (replacing sextics with their trigonal models) is
the fact that the number of points to keep track of reduces
from~$6$ to~$4$, which simplifies the computation of the braid
monodromy. As a first application, we show that all groups can be
generated by loops in a small neighborhood of (any) type~$\bE_6$
singular point of the curve.

In~\S\ref{S.perturbations}, we study perturbations of
sextics considered in~\S\S\ref{S.model} and~\ref{S.group}.
We confine ourselves to a few simple cases when the perturbed
group is easily found by simple local analysis. This gives $117$
new (compared to~\cite{degt.8a2}) sextics with abelian fundamental
group
and $8$ sextics of torus type. More complicated
perturbations are not necessary, as the resulting sextics are not
new, see Remark~\ref{not.new}.

\section{The trigonal model\label{S.model}}

\subsection{Trigonal curves}
Recall
that the \emph{Hirzebruch surface}~$\Sigma_2$ is a
geometrically ruled rational
surface with an exceptional section~$E$ of
self-intersection~$(-2)$. A \emph{trigonal curve} is a reduced
curve
$\B\subset\Sigma_2$ disjoint from~$E$ and intersecting each
generic fiber of~$\Sigma_2$ at three points. A \emph{singular
fiber} (sometimes referred to as \emph{vertical tangent})
of a trigonal curve~$\B$ is a fiber of~$\Sigma_2$ that is
not transversal to~$\B$.
The double covering~$X$ of~$\Sigma_2$
ramified at $\B+E$ is an elliptic surface, and the singular fibers
of~$\B$ are the projections of those of~$X$. For this reason, to
describe the topological types of singular fibers of~$\B$,
we use (one of)
the standard notation for the types of singular elliptic fibers,
referring to the corresponding extended Dynkin diagrams. The types
are as follows:
\Dashes
\dash
$\tilde\bA_0^*$: a simple vertical tangent;
\dash
$\tilde\bA_0^{**}$: a vertical inflection tangent;
\dash
$\tilde\bA_1^*$: a node of~$\B$ with one of the branches vertical;
\dash
$\tilde\bA_2^*$: a cusp of~$\B$ with vertical tangent;
\dash
$\tilde\bA_p$, $\tilde\bD_q$, $\tilde\bE_6$, $\tilde\bE_7$,
$\tilde\bE_8$: a simple singular point of~$\B$ of the same type
with minimal possible local intersection index with the fiber.
\endDashes
For the relation to Kodaira's classification of singular elliptic
fibers and further details and references,
see~\cite{degt.kplets}. In the
present paper, we merely use the notation.

The \emph{\rom(functional\rom) $j$-invariant}
$j=j_{\B}\:\Cp1\to\Cp1$ of
a trigonal curve $\B\subset\Sigma_2$ is defined as the analytic
continuation of the function sending a point
$b$ in the base $\Cp1$ of~$\Sigma_2$ representing a
nonsingular fiber~$F$ of~$\B$ to the $j$-invariant (divided
by~$12^3$)
of the elliptic
curve covering~$F$ and ramified at $F\cap(\B+E)$. The curve~$\B$
is called \emph{isotrivial} if $j_{\B}=\const$. Such curves can
easily be enumerated, see, \eg,~\cite{degt.kplets}. The curve~$\B$
is called \emph{maximal} if it has the following properties:
\Dashes
\dash
$\B$ has no singular fibers of type~$\bD_4$;
\dash
$j=j_{\B}$ has no critical values other than~$0$, $1$, and~$\infty$;
\dash
each point in the pull-back $j^{-1}(0)$ has ramification index at
most~$3$;
\dash
each point in the pull-back $j^{-1}(1)$ has ramification index at
most~$2$.
\endDashes
The maximality of a non-isotrivial trigonal curve $\B\subset\Sigma_2$
can easily be detected by
applying the Riemann--Hurwitz formula to the map
$j_{\B}\:\Cp1\to\Cp1$; it depends only
on the (combinatorial) set of singular fibers of~$\B$,
see~\cite{degt.kplets} for details. The classification of such
curves reduces to a combinatorial problem;
a partial classification of
maximal trigonal curves in~$\Sigma_2$ is found
in~\cite{symmetric}. An important property of maximal trigonal
curves is their rigidity, see~\cite{degt.kplets}: any small
deformation of such a curve~$\B$ is isomorphic to~$\B$. For this
reason, we do not need to keep parameters in the equations below.

\subsection{The trigonal curve~$\B$}\label{s.trigonal}
Let~$B$ be an \emph{irreducible} sextic of torus type with
\emph{simple singularities only} and with at least two
type~$\bE_6$ singular point. (Below, we show that the emphasized
properties hold automatically, see~\ref{proof.moduli}.)
Clearly,
the set of inner singularities of~$B$ can only be $(3\bE_6)$,
$(2\bE_6\splus\bA_5)$, or $(2\bE_6\splus2\bA_2)$. Hence,
according to~\cite{symmetric}, $B$
has an involutive symmetry (\ie, projective automorphism)~$c$
stable under equisingular deformations. Let~$L_c$ and~$O_c$ be,
respectively, the fixed line and the isolated fixed point of~$c$.
One has $O_c\notin B$.
Denote
by $\Cp2(O_c)$ the blow-up of~$\Cp2$ at~$O_c$. Then,
the quotient $\Cp2(O_c)/c$ is the
Hirzebruch surface~$\Sigma_2$
and the projection $B/c$ is a trigonal curve $\B\subset\Sigma_2$
with the set of singularities
$\bE_6\splus\bA_2$. The double covering $\Cp2(O_c)\to\Sigma_2$ is
ramified at~$E$ and a generic section $\L\subset\Sigma_2$ (the image
$L_c/c$) disjoint from~$E$ and not passing through the
type~$\bE_6$ singular point of~$\B$ (as otherwise the two
type~$\bE_6$ singular points of~$B$ would merge to a single
non-simple singularity).

Conversely, given a trigonal curve $\B\subset\Sigma_2$ with the
set of singularities $\bE_6\splus\bA_2$ and a section
$\L\subset\Sigma_2$ disjoint from~$E$ and not passing through the
type~$\bE_6$ singular point of~$\B$, the pull-back of~$\B$ in the
double covering of $\Sigma_2/E$ ramified at $E/E$ and~$\L$ is a
sextic $B\subset\Cp2$ with at least two type~$\bE_6$ singular
points. Below we show that $B$ is necessarily of torus type,
see~\eqref{eq.torus}.

\subsection{Equations}\label{s.equations}
Any trigonal curve $\B\subset\Sigma_2$ with the set of
singularities $\bE_6\splus\bA_2$ is either isotrivial or maximal
(see~\cite{symmetric} for precise definitions); in particular,
such curves are rigid, \ie, within each of the two families, any
two curves are isomorphic in~$\Sigma_2$. A curve~$\B$ can be
obtained by an elementary transformation from a cuspidal cubic
$C\subset\Sigma_1=\Cp2(O)$: the blow-up center~$O$ should be chosen
on the inflection tangent to~$C$, and the elementary
transformation should contract this tangent.

In appropriate affine coordinates $(x,y)$ in~$\Sigma_2$ any
trigonal
curve~$\B$ as above can be given by an equation of the form
$$
f_r(x,y):=y^3+r^2y^2+2rxy+x^2=0,
\eqtag\label{eq.equation}
$$
where $r\in\C$ is a parameter. If $r=0$, the curve is isotrivial,
its $j$-invariant being $j\equiv0$. Otherwise, the automorphism
$(x,y)\mapsto(r^3x,r^2y)$ of~$\Sigma_2$ converts the curve to
$f_1(x,y)=0$. Below, in all plots and numeric evaluation, we use
the value $r=3$.

The
$y$-discriminant
of the polynomial~$f_r$ given
by~\eqref{eq.equation} is $-x^3(27x-4r^3)$. Thus, if $r\ne0$, the
curve has three singular fibers, of types~$\tilde\bA_2$,
$\tilde\bA_0^{*}$ (vertical tangent),
and~$\tilde\bE_6$ over $x=0$, $4r^3\!/27$,
and~$\infty$, respectively. In the isotrivial case $r=0$, there
are two singular fibers, of types~$\tilde\bA_2^*$
and~$\tilde\bE_6$, over $x=0$ and~$\infty$, respectively.

The curve~$\B$ is rational; it can be parameterized by
$$
x=\xt:=rt^2+t^3,\qquad y=\yt:=-t^2.
\eqtag\label{eq.parameterization}
$$
The vertical tangency point of~$\B$ corresponds to the value
$t=-2r/3$.

Consider a section~$\L$ of~$\Sigma_2$ given by
$$
y=s(x):=ax^2+bx+c,\quad a\ne0.
\eqtag\label{eq.section}
$$
(The assumption $a\ne0$ is due to the fact that $\L$ should not
pass through the type~$\bE_6$ singular point of~$\B$.) Let
$B\subset\Cp2$ be the pull-back of~$\B$ under the double covering
of $\Sigma_2/E$ ramified at $E/E$ and~$\L$. It is a plane sextic
which, in appropriate affine coordinates $(x,y)$ in~$\Cp2$, is
given by the equation
$$
f_r(x, y^2+s(x))=0.
\eqtag\label{eq.sextic}
$$
Obviously, $B$ is of torus type, the torus structure being
$$
f_r(x,\bar y)=\bar y^3+(r\bar y+x)^2,\quad \bar y=y^2+s(x).
\eqtag\label{eq.torus}
$$
According to~\cite{degt.Oka}, this is the only torus structure
on~$B$. The inner singularities of~$B$ are two type~$\bE_6$ points
over the type~$\bE_6$ point of~$\B$ and two cusps or one
type~$\bA_5$ or~$\bE_6$ point over the cusp of~$\B$. (There is
only one point if $\L$ passes through the cusp of~$\B$; this point
is of type~$\bE_6$ if $\L$ is tangent to~$\B$ at the cusp.) The
outer singularities of~$B$ arise from the tangency of~$\L$
and~$\B$: each point of $p$-fold intersection, $p>1$, of~$\L$
and~$\B$ smooth for~$\B$ gives rise to a type~$\bA_{p-1}$ outer
singularity of~$B$. For detail, see~\cite{degt.Oka3}.

In the rest of this section, we discuss various degenerations of
the pair $(\B,\L)$ and parameterize the corresponding triples
$(a,b,c)$. For convenience, each time we mention parenthetically
the set of singularities of
the sextic~$B$ arising from $(\B,\L)$.

\subsection{Tangents and double tangents}\label{s.tangent}
Equating the values and the derivatives of $s(\xt(t))$
and~$\yt(t)$, one concludes that
a section~$\L$ as in~\eqref{eq.section}
is tangent to~$\B$ at a point $(\xt(t),\yt(t))$,
$t\ne0$, $-2r/3$, (the set of singularities
$(2\bE_6\splus2\bA_2)\splus\bA_1$)
if and only if
$$
b=-2t^2(t+r)a-\frac2{3t+2r},\quad
c=t^4(t+r)^2a-\frac{t^3}{3t+2r}.
\eqtag\label{eq.tangent}
$$

Double tangents are described by the following lemma.

\lemma\label{double.tangent}
There exists a section~$\L$ tangent to the curve~$\B$
at two distinct points
$(\xt(t_1),\yt(t_1))$ and $(\xt(t_2),\yt(t_2))$, $t_1\ne t_2$, if
and only if $t_1+t_2=-r/3$ and neither~$t_1$ nor~$t_2$ is~$0$,
$-r/6$, or $-2r/3$.
\endlemma

\proof
Substituting $t=t_1$ and $t=t_2$ to~\eqref{eq.tangent}, equating
the resulting values of~$b$ and~$c$, solving both equations
for~$a$, and equating the results, one obtains
$(t_1-t_2)^2(3t_1+3t_2+r)=0$; now, the statement is immediate.
\endproof

Thus, a section~$\L$ as in~\eqref{eq.section} is double tangent
to~$\B$ (the set of singularities
$(2\bE_6\splus2\bA_2)\splus2\bA_1$)
if and only if, for some $t\ne0$, $-r/6$, $-2r/3$, one has
$$
\gathered
a=-\frac{27}{(3t-r)^2(3t+2r)^2},\\
b=\frac{2r(27t^2+9rt-2r^2)}{(3t-r)^2(3t+2r)^2},\\
c=-\frac{2t^3(3t+r)^3}{(3t-r)^2(3t+2r)^2}.
\endgathered
\eqtag\label{eq.double}
$$

A point of quadruple intersection of~$\L$ and~$\B$ can be obtained
from Lemma~\ref{double.tangent} letting $t_1=t_2$. (Alternatively,
one can equate the derivatives of order~$0$ to~$3$ of $s(\xt(t))$
and $\yt(t)$.) As a result, $(\xt(t),\yt(t))$ is a point of
quadruple intersection of~$\L$ and~$\B$ (the set of singularities
$(2\bE_6\splus2\bA_2)\splus\bA_3$) if and only if
$$
t=-\frac{r}6,\qquad
(a,b,c)=\Bigl(-\frac{16}{3r^4},-\frac{88}{81r},\frac{r^2}{4374}\Bigr).
\eqtag\label{eq.2e6+2a2+a3}
$$
All points of intersection of this section~$\L$ and~$\B$ are:
\Dashes
\dash
transversal intersection at
$t=\Bigl(-\dfrac23+\dfrac{\sqrt2}2\Bigr)r$,
$x=\Bigl(-\dfrac{19}{54}+\dfrac{\sqrt2}4\Bigr)r^3
\approx.0459$;
\dash
transversal intersection at
$t=\Bigl(-\dfrac23-\dfrac{\sqrt2}2\Bigr)r$,
$x=\Bigl(-\dfrac{19}{54}-\dfrac{\sqrt2}4\Bigr)r^3
\approx-19.1$;
\dash
quadruple intersection at $t=-\dfrac{r}6$,
$x=\dfrac{5r^3}{216}=.625$.
\endDashes
The curve~$\B$ and the section~$\L$ given by~\eqref{eq.2e6+2a2+a3}
are plotted in Figure~\ref{fig.2e6+2a2+a3} (in black and grey,
respectively). The section is above the curve over $x=0$; it
intersects the topmost branch over $x\approx.0459$ and is tangent
to the middle branch over $x=.625$.

\midinsert
\plot{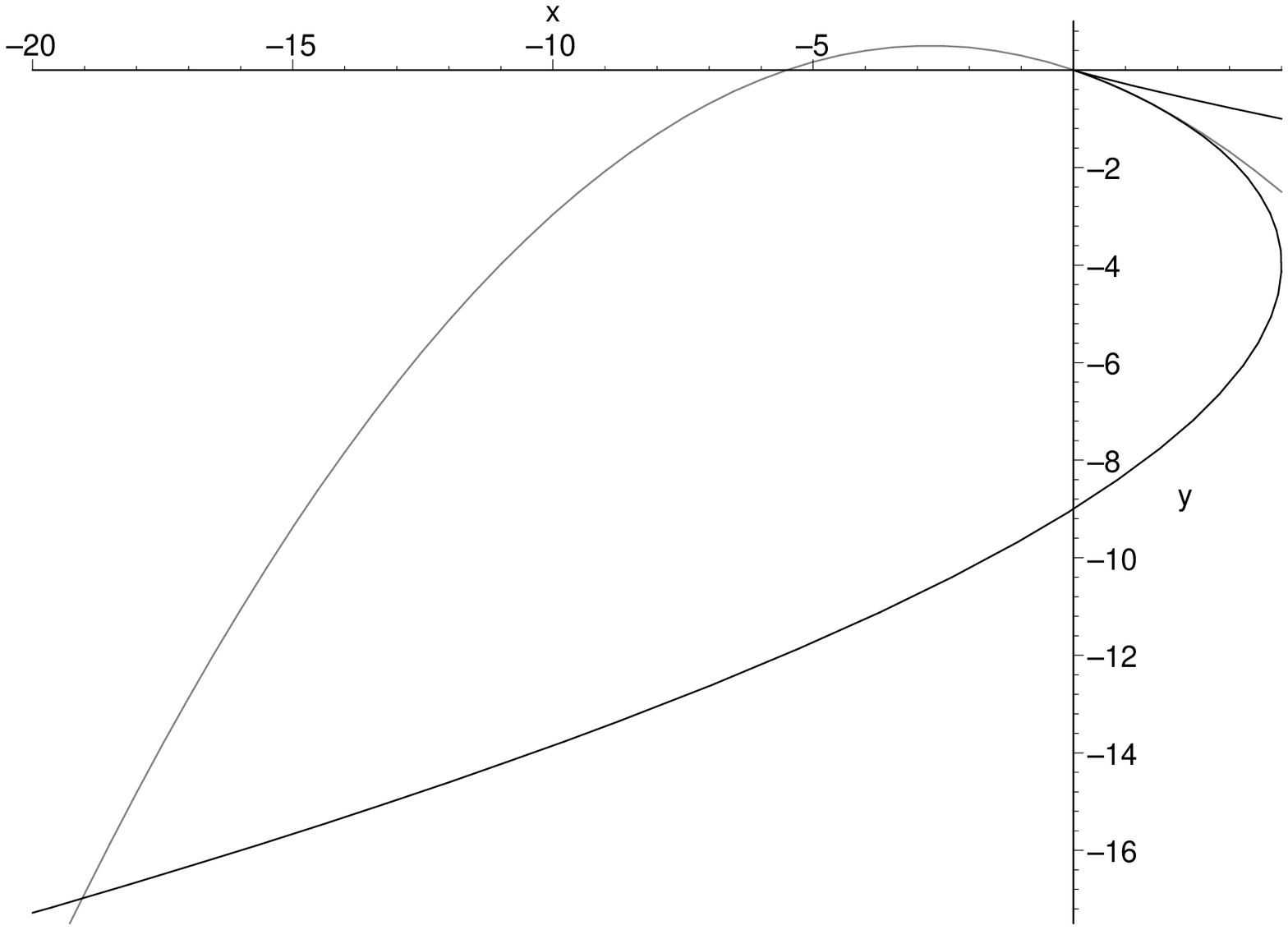}
\figure
The set of singularities $(2\bE_6\splus2\bA_2)\splus\bA_3$
\endfigure\label{fig.2e6+2a2+a3}
\endinsert

\subsection{Sections through the cusp}\label{s.cusp}
A section~$\L$ as in~\eqref{eq.section}
passes through the cusp of~$\B$ (the set of
singularities $(2\bE_6\splus\bA_5)$)
if and only if $c=0$;
it is tangent to~$\B$ at the cusp (the set of singularities
$(3\bE_6)$)
if and only if, in addition,
$b=-1/r$.

A section tangent to~$\B$ at a point $(\xt(t),\yt(t))$,
see~\eqref{eq.tangent},
passes through the cusp of~$\B$ (the set of
singularities $(2\bE_6\splus\bA_5)\splus\bA_1$)
if and only if
$$
a=\frac1{t(t+r)^2(3t+2r)},\quad
b=-\frac{2(2t+r)}{(t+r)(3t+2r)},\quad
c=0,
\eqtag\label{eq.cusp.tangent}
$$
$t\ne0$, $-r$, $-2r/3$. (Note that the value $t=-r$ corresponds to
the smooth point of~$\B$ in the same vertical fiber as the cusp.)
Such a section is tangent to~$\B$ at the cusp (the set of
singularities $(3\bE_6)\splus\bA_1$) if and only if
$$
t=-\frac{r}3,\qquad
(a,b,c)=\Bigl(-\frac{27}{4r^4},-\frac1r,0\Bigr).
\eqtag\label{eq.3e6+a1}
$$
The points of intersection of the latter section~$\L$ and~$\B$
are:
\Dashes
\dash
the cusp of~$\B$ at $t=0$, $x=0$;
\dash
transversal intersection at $t=-\dfrac{4r}3$,
$x=-\dfrac{16r^3}{27}=-16$;
\dash
tangency at $t=-\dfrac{r}3$, $x=\dfrac{2r^3}{27}=2$.
\endDashes
The section~$\L$ given by~\eqref{eq.3e6+a1} looks similar to that
shown in Figure~\ref{fig.2e6+2a2+a3}. (Near the cusp of~$\B$, the
two curves are too close to be distinguished visually.) Between
$x=0$ and $x=2$, the section lies between the topmost and middle
branches of~$\B$.

\subsection{Inflection tangents}\label{s.inflection}
Equating the derivatives of $s(\xt(t))$ and $\yt(t)$ of order $0$,
$1$, and~$2$, one can see that
a section~$\L$ as in~\eqref{eq.section}
is inflection tangent to~$\B$ at a
point $(\xt(t),\yt(t))$, $t\ne0$, $-2r/3$, (the set of
singularities $(2\bE_6\splus2\bA_2)\splus\bA_2$) if and only if
$$
a=\frac3{t(3t+2r)^3},\quad
b=-\frac{2(12t^2+15rt+4r^2)}{(3t+2r)^3},\quad
c=-\frac{t^3(6t^2+6rt+r^2)}{(3t+2r)^3}.
\eqtag\label{eq.inflection}
$$
Such a section passes through the cusp of~$\B$ (the set of
singularities $(2\bE_6\splus\bA_5)\splus\bA_2$) if and only if
$t=(-3\pm\sqrt3)r/6$. Thus, we obtain \emph{two families}, which
are Galois conjugate over $\Q[\sqrt3]$, \cf.~\cite{OkaPho.moduli}.
For one of the families, one has
$$
t=\Bigl(-\frac12+\frac{\sqrt3}6\Bigr)r,\quad
(a,b,c)=\Bigl(\frac{12(3-2\sqrt3)}{r^4},-\frac{4(2-\sqrt3)}r,0\Bigr),
\eqtag\label{eq.2e6+a5+a2.1}
$$
and the points of intersection of~$\L$ and~$\B$ are:
\Dashes
\dash
the cusp of~$\B$ at $t=0$, $x=0$;
\dash
transversal intersection at
$t=\Bigl(-\dfrac12-\dfrac{\sqrt3}2\Bigr)r$,
$x=\Bigl(-\dfrac14-\dfrac{\sqrt3}4\Bigr)r^3\approx-18.4$;
\dash
inflection tangency at $t=\Bigl(-\dfrac12+\dfrac{\sqrt3}6\Bigr)r$,
$x=\Bigl(\dfrac1{12}-\dfrac{\sqrt3}{36}\Bigr)r^3\approx.951$.
\endDashes
This section looks similar to that shown in
Figure~\ref{fig.2e6+2a2+a3}; between $x=0$ and $x\approx.951$, the
section is just below the middle branch of the curve.

For the other family, one has
$$
t=\Bigl(-\frac12-\frac{\sqrt3}6\Bigr)r,\quad
(a,b,c)=\Bigl(\frac{12(3+2\sqrt3)}{r^4},-\frac{4(2+\sqrt3)}r,0\Bigr),
\eqtag\label{eq.2e6+a5+a2.2}
$$
and the points of intersection of~$\L$ and~$\B$ are:
\Dashes
\dash
the cusp of~$\B$ at $t=0$, $x=0$;
\dash
transversal intersection at
$t=\Bigl(-\dfrac12+\dfrac{\sqrt3}2\Bigr)r$,
$x=\Bigl(-\dfrac14+\dfrac{\sqrt3}4\Bigr)r^3\approx4.94$;
\dash
inflection tangency at $t=\Bigl(-\dfrac12-\dfrac{\sqrt3}6\Bigr)r$,
$x=\Bigl(\dfrac1{12}+\dfrac{\sqrt3}{36}\Bigr)r^3\approx3.55$.
\endDashes
The curve~$\B$ and the section~$\L$ given by~\eqref{eq.2e6+a5+a2.2}
are plotted in Figure~\ref{fig.2e6+a5+a2}, in black and grey,
respectively.

\midinsert
\plot{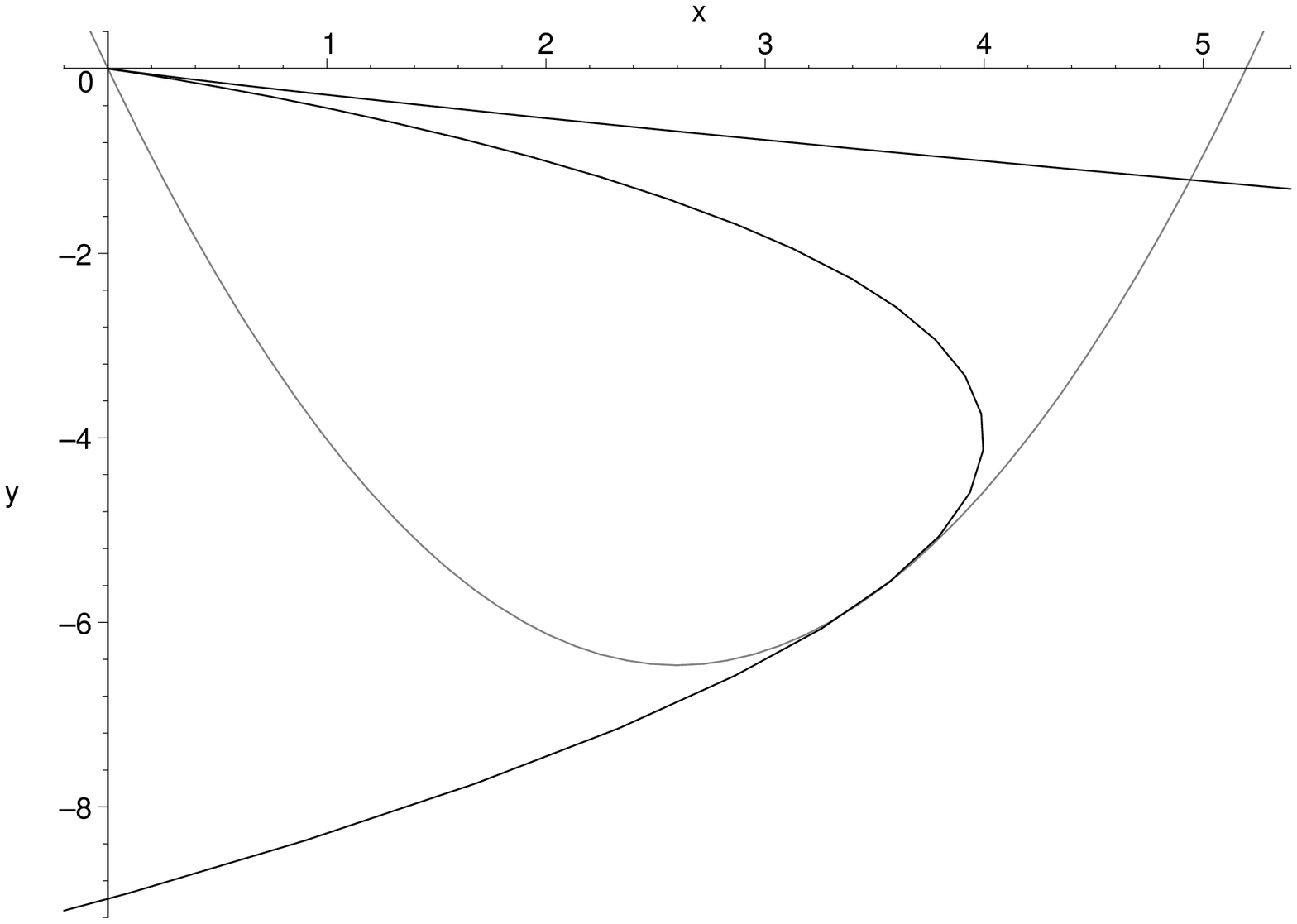}
\figure
The set of singularities
$(2\bE_6\splus\bA_5)\splus\bA_2$, the family~\eqref{eq.2e6+a5+a2.2}
\endfigure\label{fig.2e6+a5+a2}
\endinsert

\subsection{Proof of Theorem~\ref{th.moduli}}\label{proof.moduli}
First, note that any
sextic with two type~$\bE_6$ singular points
is irreducible and has simple singularities only. The first
statement follows from the fact that an irreducible curve of
degree~$4$ or~$5$ (respectively,~$\le3$) may have at most one
(respectively, none) type~$\bE_6$ singular point, and the second
one, from the fact that a type~$\bE_6$ (respectively, non-simple)
singular point takes~$3$ (respectively,~$\ge6$) off the
genus, whereas the genus of a nonsingular sextic is~$10$. Thus, we
can apply the results of~\cite{symmetric} enumerating stable
symmetries of curves.

For a set of
singularities~$\Sigma\supset2\bE_6$, consider the moduli
space $\CM(\Sigma)$ of sextics~$B$ of torus type with the set of
singularities~$\Sigma$ and
the moduli space $\tCM(\Sigma)$ of pairs $(B,c)$, where $B$
is a sextic as above and $c$ is a
stable involution of~$B$. Due to~\cite{symmetric}, the forgetful
map $\tCM(\Sigma)\to\CM(\Sigma)$ is generically finite-to-one and
onto.

As explained in Sections~\ref{s.trigonal} and~\ref{s.equations},
the space $\tCM(\Sigma)$ can be identified with the moduli space
of pairs $(\B,\L)$, where $\B\subset\Sigma_2$ is a trigonal curve
given by~\eqref{eq.equation} and $\L$ is a section of~$\Sigma_2$
in a certain prescribed position with respect to~$\B$. The spaces
of pairs $(\B,\L)$ are described in
Sections~\ref{s.tangent}--\ref{s.inflection}, and for each
$\Sigma\ne(2\bE_6\splus\bA_5)\splus\bA_2$, an explicit rational
parameterization is found. (Strictly speaking, in order to pass to
the moduli, we need to fix a value of~$r$, say, $r=3$. This
results in a Zariski open subset of the moduli space. The portion
corresponding to $r=0$ has positive codimension as the isotrivial
curve $f_0=0$ has $1$-dimensional group~$\C^*$ of symmetries.)
Hence, the space $\tCM(\Sigma)$ is rational and, if
$\dim\CM(\Sigma)\le2$, so is $\CM(\Sigma)$. The only case when
$\dim\CM(\Sigma)\ge3$ is $\Sigma=(2\bE_6\splus2\bA_2)$. In this
case, each curve~$B$ has a unique stable involution,
see~\cite{symmetric}, and the map $\tCM(\Sigma)\to\CM(\Sigma)$ is
generically one-to-one; hence, $\CM(\Sigma)$ is still rational.

In the exceptional case $\Sigma=(2\bE_6\splus\bA_5)\splus\bA_2$, the
space $\CM(\Sigma)=\tCM(\Sigma)$ consists of two points. The fact
that any sextic with this set of singularities is of torus type
follows immediately from~\cite{JAG}.
\qed

\Remark
The only sets of singularities containing $2\bE_6$
where the curves have more than one
(three) stable involutions are $(3\bE_6)$ and $(3\bE_6)\splus\bA_1$,
see~\cite{symmetric}. In both cases, the group of stable
symmetries can be identified with the group~$\SG3$ of permutations
of the three type~$\bE_6$ points. It follows that all three
involutions are conjugate by stable symmetries; hence, the map
$\tCM(\Sigma)\to\CM(\Sigma)$ is still one-to-one.
\endRemark

\section{Van Kampen's method in Hirzebruch surfaces}\label{S.vanKampen}

In this section, we give a formal and
detailed exposition of a few
observations outlined in~\cite{degt.Oka3}. Keeping in mind future
applications, we treat the general case of a Hirzebruch
surface~$\Sigma_k$, $k\ge1$, and a \emph{$d$-gonal curve}
$C\subset\Sigma_k$, see Definition~\ref{def.d-gonal}.

Certainly, the essence of this approach is due to van
Kampen~\cite{vanKampen}; we merely introduce a few restrictions to
the objects used in the construction which make the choices
involved slightly more canonical and easier to handle. By no means
do we assert that the restrictions are necessary for the approach
to work in general.

\subsection{Preliminary definitions}\label{s.defs}
Fix a Hirzebruch surface~$\Sigma_k$, $k\ge1$. Denote by
$p\:\Sigma_k\to\Cp1$ the ruling, and let $E\subset\Sigma_k$ be
the exceptional section, $E^2=-k$. Given a point~$b$ in the
base~$\Cp1$, we denote by~$F_b$ the fiber $p^{-1}(b)$. Let
$F_b^\circ$ be the `open fiber' $F_b\sminus E$. Observe that
$F_b^\circ$ is a dimension~$1$ affine space over~$\C$; hence, one
can speak about lines, circles, convexity, convex hulls, \etc.
in~$F_b^\circ$.
(Thus, strictly speaking, the notation~$F_b^\circ$ means slightly
more than just the set theoretical difference $F_b\sminus E$: we
always consider~$F_b^\circ$ with its canonical affine structure.)
Define the \emph{convex hull} $\Conv C$
of a subset
$C\subset\Sigma_k\sminus E$ as the union of its fiberwise convex
hulls:
$$
\Conv C=\bigcup_{b\in\Cp1}\Conv(C\cap F_b^\circ).
$$

\definition\label{def.d-gonal}
Let $d\ge1$ be an integer.
A \emph{$d$-gonal curve} (or \emph{degree~$d$ curve})
on~$\Sigma_k$ is a reduced algebraic
curve $C\in\ls|dE+dkF|$ disjoint from
the exceptional section~$E$. (Here, $F$ is any fiber
of~$\Sigma_k$.) A \emph{singular fiber} of a $d$-gonal curve~$C$
is a fiber of~$\Sigma_k$ that intersects~$C$ at fewer than $d$
points. (With a certain abuse of the language, the points in the
base $\Cp1$ whose pull-backs are singular fibers will also be
referred to as singular fibers of~$C$.)
\enddefinition

\Remark\label{affine}
Recall that the complement $\Sigma_x\sminus E$ can be covered by
two affine charts, with coordinates $(x,y)$ and $(x',y')$ and
transition function $x'=1/x$, $y'=y/x^k$. In the coordinates
$(x,y)$, any $d$-gonal curve~$C$ is given by an equation of the form
$$
f(x,y)=\sum_{i=0}^dy^iq_{i}(x)=0,\quad
 \deg q_i=k(d-i),\quad
 q_d=\const\ne0,
$$
and the singular fibers of~$C$ are those of the form $F_x$, where $x$ is
a root of the $y$-discriminant~$D_y$ of~$f$. (The fiber $F_\infty$ over
$x=\infty$ is singular for~$C$ if and only if $\deg D_y<kd(d-1)$.)
\endRemark

\subsection{Proper sections and braid monodromy}\label{s.proper}
Fix a $d$-gonal curve $C\subset\Sigma_k$.

\definition\label{def.proper}
Let $\disk\subset\Cp1$ be a closed (topological) disk. A partial
section $s\:\disk\to\Sigma_k$ of~$p$ is called \emph{proper} if its
image is disjoint from both~$E$ and $\Conv C$.
\enddefinition

\lemma\label{proper.homotopic}
Any disk $\disk\subset\Cp1$ admits a proper section
$s\:\disk\to\Sigma_k$.
Any two proper sections over~$\disk$ are homotopic
in the class of proper sections\rom;
furthermore, any homotopy over a fixed
point $b\in\disk$ extends to a homotopy over~$\disk$.
\endlemma

\proof
The restriction~$p'$ of~$p$ to $\Sigma_k\sminus(E\cup\Conv C)$ is a
locally trivial fibration with a typical fiber~$F'$
homeomorphic to a punctured open disk. Since $\disk$ is
contractible, $p'$ is trivial over~$\disk$ and, after trivializing,
sections over~$\disk$ can be identified with maps $\disk\to F'$.
Such maps do exist, and
any two such maps are homotopic, again due to the fact that
$\disk$ is contractible.
\endproof

Pick a closed disk $\disk\subset\Cp1$ as above and denote
$\disk^\sharp=\disk\sminus\{b_1,\ldots,b_l\}$,
where $b_1,\ldots,b_l$ are the singular fibers of~$C$ that belong
to~$\disk$.
Fix a point
$b\in\disk^\sharp$. The restriction
$p^\sharp\:p^{-1}(\disk^\sharp)\sminus(C\cup E)\to\disk^\sharp$
is a locally trivial fibration with a typical fiber
$F_b^\circ\sminus C$, and any proper section $s\:\disk\to\Sigma_k$
restricts to
a section of~$p^\sharp$. Hence, given a proper
section~$s$, one can define the group
$\pi_F:=\pi_1(F_b^\circ\sminus C,s(b))$ and the \emph{braid
monodromy} $m\:\pi_1(\disk^\sharp,b)\to\Aut\pi_F$. Informally, for
a loop $\Gs\:[0,1]\to\disk^\sharp$, the automorphism $m([\Gs])$
of~$\pi_F$ is
obtained by dragging the fiber~$F_b$ along $\Gs(t)$ while keeping the
base point on $s(\Gs(t))$. (Formally, it is obtained by
trivializing the fibration $\Gs^*p^\sharp$.)

It is essential that, in this paper,
we reserve the term
`braid monodromy' for the homomorphism~$m$ constructed using a
\emph{proper} section~$s$.
Under this convention,
the following
lemma is an immediate consequence of
Lemma~\ref{proper.homotopic} and the obvious fact that the braid
monodromy is homotopy invariant.

\lemma\label{proper.braid}
The braid monodromy $m\:\pi_1(\disk^\sharp,b)\to\Aut\pi_F$ is well
defined and independent of the choice of a proper section
over~$\disk$ passing through $s(b)$.
\qed
\endlemma

\Remark
More generally, given a path
$\tilde\Gs\:[0,1]\to p^{-1}(\disk^\sharp)\sminus(\Conv C\cup E)$,
one can use Lemma~\ref{proper.homotopic} to conclude that the
braid monodromy commutes with the translation isomorphism
$$
T_{\Gs}\:\pi_1(\disk^\sharp,\Gs(0))\to\pi_1(\disk^\sharp,\Gs(1))
$$
(where $\Gs=p\circ\tilde\Gs\:[0,1]\to\disk^\sharp$)
and the isomorphism
$$
\Aut\pi_1(F_{\Gs(0)}^\circ\sminus C,\tilde\Gs(0))\to
 \Aut\pi_1(F_{\Gs(1)}^\circ\sminus C,\tilde\Gs(1))
$$
induced by the translation $T_{\tilde\Gs}$ along~$\tilde\Gs$.
\endRemark

\Remark\label{constant.section}
For most computations, we will take for~$s$ a `constant section'
constructed as follows: pick an affine coordinate system $(x,y)$,
see Remark~\ref{affine}, so that the point $x=\infty$ does
\emph{not} belong to~$\disk$, and let~$s$ be the section
$x\mapsto c=\const$, $\ls|c|\gg0$. (In other words, the graph
of~$s$ is the $1$-gonal curve $\{y=c\}\subset\Sigma_k$.) Since the
intersection $p^{-1}(\disk)\cap\Conv C\subset\Sigma_k\sminus E$ is
compact, such a section is indeed proper whenever $\ls|c|$ is
sufficiently large.
\endRemark

\Remark\label{extension}
Another consequence of Lemma~\ref{proper.braid} is the fact
that, for
any nested pair of disks $\disk_1\subset\disk_2$, the braid
monodromy commutes with the inclusion homomorphism
$\pi_1(\disk_1^\sharp)\to\pi_1(\disk_2^\sharp)$. Indeed, one can
construct both
monodromies using a proper section
over~$\disk_2$ and restricting it to~$\disk_1$ when necessary.
\endRemark

Pick a basis $\zeta_1,\ldots,\zeta_d$ for $\pi_F$ and a basis
$\Gs_1,\ldots,\Gs_l$ for $\pi_1(\disk^\sharp,b)$. Denote
$m_i=m(\Gs_i)$, $i=1,\ldots,l$.
The following
statement is the essence of Zariski--van Kampen's method for
computing the fundamental group of a plane algebraic curve,
see~\cite{vanKampen} for the proof and further details.

\theorem\label{th.vanKampen.original}
Let $\disk\subset\Cp1$ be a closed disk as above,
and assume that the boundary
$\partial\disk$ is free of singular fibers of~$C$. Then one has
$$
\pi_1(p^{-1}(\disk)\sminus(C\cup E),s(b))=
 \bigl<\zeta_1,\ldots,\zeta_d\bigm|m_i=\id,\ i=1,\ldots,l\bigl>,
$$
where each \emph{braid relation} $m_i=\id$ should be understood as
a $d$-tuple of relations $\zeta_j=m_i(\zeta_j)$, $j=1,\ldots,d$.
\qed
\endtheorem

\subsection{The monodromy at infinity}\label{s.monodromy.infty}
Let $b\in\disk^\sharp\subset\disk\subset\Cp1$ be as in
Section~\ref{s.proper}. Denote by $\Gr_b\in\pi_F$ the
`counterclockwise' generator of the abelian subgroup
$\Z\cong\pi_1(F_b^\circ\sminus\Conv C)\subset\pi_F$. (In other
words, $\Gr_b$ is the class of a large circle in~$F_b^\circ$
encompassing $\Conv C\cap F_b^\circ$.
If
$\zeta_1,\ldots,\zeta_d$ is a `standard basis' for~$\pi_F$, \cf.
Figure~\ref{fig.basis}, left, then
$\Gr_b=\zeta_1\cdot\ldots\cdot\zeta_d$.) Clearly, $\Gr_b$ is invariant
under the braid monodromy and, properly understood, it is
preserved by the translation homomorphism along any path in
$p^{-1}(\disk^\sharp)\sminus(\Conv C\cup E)$. (Indeed,
as explained in the proof of Lemma~\ref{proper.homotopic},
the fibration
$p^{-1}(\disk)\sminus(\Conv C\cup E)\to\disk$ is trivial,
hence $1$-simple.) Thus, there is a canonical identification of
the elements $\Gr_{b'}$, $\Gr_{b''}$ in the fibers over any two
points $b',b''\in\disk^\sharp$; for this reason, we will omit the
subscript~$b$ in the sequel.

Assume that the boundary $\partial\disk$ is free of singular
fibers of~$C$. Then, connecting $\partial\disk$ with the base
point~$b$ by a path in~$\disk^\sharp$ and traversing it in the
counterclockwise direction (with respect to the canonical complex
orientation of~$\disk$), one obtains a certain element
$[\partial\disk]\in\pi_1(\disk^\sharp,b)$ (which depends
on the choice of the path above).

\proposition\label{monodromy.infty}
In the notation above, assume that the interior of~$\disk$
contains all singular fibers of~$C$. Then, for any
$\zeta\in\pi_F$, one has
$m([\partial\disk])(\zeta)=\Gr^{k}\zeta\Gr^{-k}$. \rom(In
particular, $m([\partial\disk])$ does not depend on the choices in
the definition of $[\partial\disk]$.\rom)
\endproposition

\proof
Due to the homotopy invariance of the braid monodromy (and the
invariance of~$\Gr$), one can replace~$\disk$ with any larger disk
and assume that the base point~$b$ is in the boundary.
Consider affine charts $(x,y)$ and $(x',y')$, see
Remark~\ref{affine},
such that the fiber
$\{x=\infty\}=\{x'=0\}$ does not belong to~$\disk$ (and hence is
nonsingular for~$C$), and replace $\disk$ with the disk
$\{\ls|x|\le1/\epsilon\}$ for some positive $\epsilon\ll1$.
About $x'=0$, the curve~$C$ has $d$ analytic branches of the form
$y'=c_i+x'\Gf_i(x')$,
where $c_i$ are pairwise distinct constants and $\Gf_i$ are
analytic functions, $i=1,\ldots,d$. Restricting these expressions
to the circle $x'=\epsilon\exp(-2\pi t)$, $t\in[0,1]$,
and passing to
$x=1/x'$ and $y=y'x^k$, one obtains
$y=c_i\epsilon^{-k}\exp(2k\pi t)+O(\epsilon^{-k+1})$,
$i=1,\ldots,d$. Thus,
from the point of view of a trivialization of the ruling
over~$\disk$ (\eg, the one given by~$y$),
the parameter~$\epsilon$ can be chosen so small that
the $d$
branches move along $d$ pairwise disjoint concentric circles
(not quite round), each
branch making $k$ turns in the counterclockwise direction. On the
other hand, one can assume that
the base point remains in a constant section $y=c=\const$
with $\ls|c|\gg\epsilon^{-k}\max\ls|c_i|$, see
Remark~\ref{constant.section}. The resulting braid is the
conjugation by~$\Gr^{-k}$.
\endproof

\subsection{The relation at infinity}
We are ready to state the principal result of this section.
Fix a $d$-gonal curve $C\subset\Sigma_k$ and
choose a closed disk $\disk\subset\Cp1$ satisfying the following
conditions:
\roster
\item\local{all.but.1}
$\disk$ contains all but at most one singular fibers of~$C$;
\item\local{no.boundary}
none of the singular fibers of~$C$ is in the boundary
$\partial\disk$.
\endroster
As in Section~\ref{s.proper}, pick a base point
$b\in\disk^\sharp$, a basis $\zeta_1,\ldots,\zeta_d$ for the group
$\pi_F$ over~$b$, and a basis $\Gs_1,\ldots,\Gs_l$ for the group
$\pi_1(\disk^\sharp,b)$. Let $m_i=m(\Gs_i)$, $i=1,\ldots,l$, where
$m\:\pi_1(\disk^\sharp,b)\to\Aut\pi_F$ is the braid monodromy.

\theorem\label{th.vanKampen}
Under the assumptions~\loccit{all.but.1}, \loccit{no.boundary}
above, one has
$$
\pi_1(\Sigma_k\sminus(C\cup E))=
 \bigl<\zeta_1,\ldots,\zeta_d\bigm|
 m_i=\id,\ i=1,\ldots,l,\ \ \Gr^k=1\bigr>,
$$
where each \emph{braid relation} $m_i=\id$ should be understood as
a $d$-tuple of relations $\zeta_j=m_i(\zeta_j)$, $j=1,\ldots,d$,
and $\Gr\in\pi_F$ is the element introduced in
Section~\ref{s.monodromy.infty}.
\endtheorem

The relation $\Gr^k=1$ in Theorem~\ref{th.vanKampen} is called the
\emph{relation at infinity}. If $k=1$, it coincides with the well
known relation $\Gr=1$ for the group of a plane curve.

\proof
First, consider the case when $\disk$ contains \emph{all} singular
fibers of~$C$. As in the proof of
Proposition~\ref{monodromy.infty}, one can replace~$\disk$ with
any larger disk, \eg, with the one given by
$\{\ls|x|\le1/\epsilon\}$, where $(x,y)$ are affine coordinates
such that the point $x=\infty$ is not in~$\disk$ and $\epsilon$ is
a sufficiently small positive real number. Furthermore, one can
take for~$s$ a constant section $x\mapsto\epsilon^{-k}c=\const$,
$\ls|c|\gg0$,
see Remark~\ref{constant.section}, and choose the base point~$b$
in the boundary $\partial\disk$.
The fundamental group
$\pi_1(p^{-1}(\disk)\sminus(C\cup E))$ is given by
Theorem~\ref{th.vanKampen.original}, and the patching of the
nonsingular fiber $\{x=\infty\}=\{x'=0\}$
results in the additional relation $[\partial\Gamma]=1$, where
$\Gamma$ is the disk $\{y'=c,\ \ls|x'|\le\epsilon\}$. (Here,
$x'=1/x$ and $y'=y/x^k$ are the affine coordinates in the
complementary chart, see Remark~\ref{affine}.
We assume that the constant $\ls|c|$ is so large that
$\Gamma\cap\Conv C=\varnothing$.) Restricting to the
boundary $x'=\epsilon\exp(-2\pi t)$, $t\in[0,1]$, and passing back
to $(x,y)$, one finds that the loop $\partial\Gamma$ is given by
$x=\epsilon^{-1}\exp(2\pi t)$, $y=\epsilon^{-k}c\exp(2k\pi t)$;
it is homotopic to $\Gr^k\cdot[s(\partial\disk)]$. Since the loop
$s(\partial\disk)$ is contractible (along the image of~$s$), the
extra relation is $\Gr^k=1$, as stated.

Now, assume that one singular fiber of~$C$ is not in~$\disk$. Extend
$\disk$ to a larger disk $\disk'\supset\disk$ containing the
missing singular fiber (and extend the braid monodromy, see
Remark~\ref{extension}). For~$\disk'$, the theorem has already
been proved, and the resulting presentation of the group differs
from the one given by~$\disk$ by an extra relation $m_{l+1}=\id$.
However, under an appropriate choice of the
additional generator~$\Gs_{l+1}$,
one has
$[\partial\disk']=[\partial\disk]\cdot\Gs_{l+1}$. Clearly,
$m([\partial\disk])$ is a word in $m_1,\ldots,m_l$ and, in view of
Proposition~\ref{monodromy.infty}, the monodromy
$m([\partial\disk'])$ is the conjugation by~$\Gr^{-k}$. Hence, in
the presence of the relation at infinity $\Gr^k=1$, the
additional relation
$m_{l+1}=\id$ is a consequence of the other braid relations,
and the statement follows.
\endproof

\section{The fundamental group\label{S.group}}

\subsection{Preliminaries}
Fix a sextic~$B$, pick a stable involutive symmetry~$c$ of~$B$,
see~\S\ref{S.model}, and let $\B,\L\subset\Sigma_2=\Cp2(O_c)/c$ be
the projections of~$B$ and~$L_c$, respectively. We start with
applying Theorem~\ref{th.vanKampen} to the $4$-gonal curve $\B+\L$
and computing the group
$\bpi:=\pi_1(\Sigma_2\sminus(\B\cup\L\cup E))$.

In order to visualize the braid monodromy, we
will consider the standard \emph{real structure} (\ie,
anti-holomorphic involution) $\conj\:(x,y)\mapsto(\bar x,\bar y)$
on~$\Sigma_2$, where bar stands for the complex conjugation. A
reduced algebraic
curve $C\subset\Sigma_2$ is said to be \emph{real} (with respect
to~$\conj$)
if it is $\conj$-invariant (as a set).
Alternatively, $C$ is real if
and only if, in
the coordinates $(x,y)$, it can be given by a polynomial with real
coefficients. In particular, the curve~$\B$ given
by~\eqref{eq.equation} is real. Given a real curve
$C\subset\Sigma_2$, one can speak about its \emph{real part}
$C_{\R}$ (\ie, the set of points of~$C$ fixed by~$\conj$), which
is a codimension~$1$ subset in the real part
of~$\Sigma_2$.

To use Theorem~\ref{th.vanKampen}, we take for~$\disk$ a closed
regular neighborhood of the smallest segment of the real axis
$\Rp1$ containing all singular fibers of $\B+\L$ except the one of
type~$\tilde\bE_6$ at infinity, see the shaded area in
Figure~\ref{fig.basis}, right.
Recall that
singular are the fiber $\{x=0\}$ through the cusp, the vertical
tangent $\{x=4\}$, and the fibers through the points of
intersection of~$\B$ and~$\L$. (As in~\S\ref{S.model}, we use the
value $r=3$ for the numeric evaluation.)
We only consider the
four extremal sections~$\L$ given by~\eqref{eq.2e6+2a2+a3},
\eqref{eq.3e6+a1}, \eqref{eq.2e6+a5+a2.1},
and~\eqref{eq.2e6+a5+a2.2}. In each case, all
singular fibers are real; they are listed in~\S\ref{S.model}.

To compute the braid monodromy, we use a constant real section
$s\:\disk\to\Sigma_2$
given by $x\mapsto\const\gg0$, see Remark~\ref{constant.section}, and
the base point $b=(\epsilon,0)\in\disk$, where $\epsilon>0$ is
sufficiently small.
The basis $\Gs_1,\ldots,\Gs_l$ for the group
$\pi_1(\disk^\sharp,b)$ is chosen as shown in
Figure~\ref{fig.basis}, right: each~$\Gs_i$ is a small loop about
a singular fiber connected to~$b$ by a real segment,
circumventing the
interfering singular fibers in the counterclockwise direction. Let
$F=F_b$ be the base fiber, and choose a basis $\Ga$, $\Gb$, $\Gg$,
$\Gd$ for the group $\pi_F=\pi_1(F^\circ\sminus(\B\cup\L),s(b))$
as shown in Figure~\ref{fig.basis}, left. (Note that, in all cases
considered below, all points of the intersection $F\cap(\B\cup\L)$
are real.)
The following notation
convention is important for the sequel.

\Remark\label{notation}
We use a double notation for the elements of the basis
for~$\pi_F$. On the one hand, to be consistent with
Theorem~\ref{th.vanKampen}, we denote them
$\zeta_1,\ldots,\zeta_4$, numbering the loops consecutively according
to the decreasing of the $y$-coordinate of the point. Then the
element $\Gr\in\pi_F$ introduced in
Section~\ref{s.monodromy.infty} is given by
$\Gr=\zeta_1\zeta_2\zeta_3\zeta_4$, and the relation at infinity in
Theorem~\ref{th.vanKampen} turns to
$(\zeta_1\zeta_2\zeta_3\zeta_4)^2=1$. On the other hand, to make
the formulas more readable, we denote the basis elements by $\Ga$,
$\Gb$, $\Gg$, and~$\Gd$. The first three elements are numbered
consecutively, whereas~$\Gd$ plays a
very special r\^ole in the passage
to the group $\pi_1(\Cp2\sminus B)$, see
Lemma~\ref{bpi->pi} below: we always assume that $\Gd$ is
the element represented
by a loop about the point $F\cap\L$.
Thus, the position of~$\Gd$ in the sequence
$(\Ga,\Gb,\Gg,\Gd)$ may change; this position is important for the
expression for~$\Gr$ and hence for the
relation at infinity.
\endRemark


\midinsert
\centerline{\picture{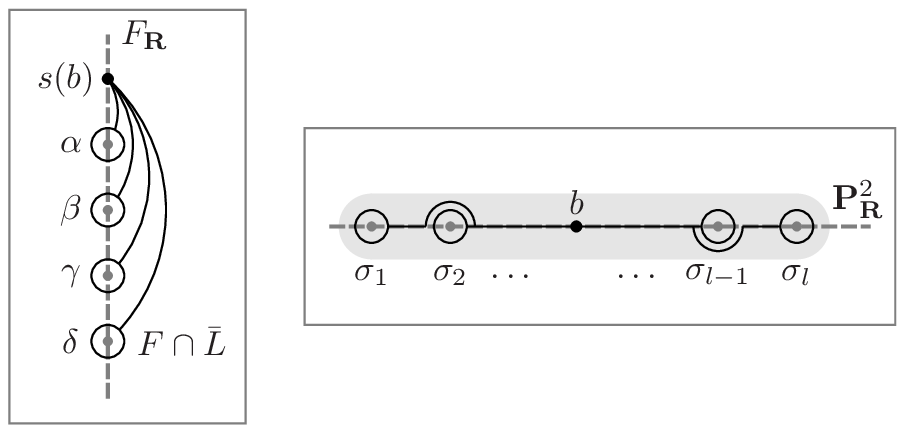}}
\figure\label{fig.basis}
The basis $\Ga$, $\Gb$, $\Gg$, $\Gd$ and the loops~$\Gs_i$
\endfigure
\endinsert

The passage from a presentation of~$\bpi$ to the
that of the group $\pi_1:=\pi_1(\Cp2\sminus B)$ is
given by the following lemma.

\lemma\label{bpi->pi}
If $\bpi$ is given by
$\bigl<\Ga,\Gb,\Gg,\Gd\bigm|R_j=1,\ j=1,\ldots,s\bigr>$, then
$$
\pi_1=\bigl<\Ga,\bGa,\Gb,\bGb,\Gg,\bGg\bigm|R'_j=\bar R'_j=1,\
 j=1,\ldots,s\bigr>,
$$
where bar stands for the conjugation by~$\Gd$,
$\bar w=\Gd^{-1}w\Gd$,
each relation $R'_j$ is obtained from $R_j$, $j=1,\ldots,s$,
by letting
$\Gd^2=1$ and expressing the result in terms of the generators
$\Ga$, $\bGa$, \dots, and $\bar R'_j=\Gd^{-1}R'_j\Gd$,
$j=1,\ldots,s$.
\rom(In other words,
$\bar R'_j$ is obtained from $R'_j$ by interchanging
$\Ga\leftrightarrow\bGa$, $\Gb\leftrightarrow\bGb$, and
$\Gg\leftrightarrow\bGg$.\rom)
\endlemma

\proof
The projection
$\Cp2\sminus(B\cup O_c)\to\Sigma_2\sminus(\B\cup E)$ is a double
covering ramified at~$\L$. Hence, one has
$\pi_1=\pi_1(\Cp2\sminus(B\cup O_c))=\Ker[\kappa\:\bpi/\Gd^2\to\CG2]$, where
$\kappa\:\Ga,\Gb,\Gg\mapsto0$ and $\kappa\:\Gd\mapsto1$. (Note
that the compactification of the
double covering above is \emph{not} ramified at~$\B$.) Lift
$\kappa$ to a homomorphism
$\tilde\kappa\:\<\Ga,\Gb,\Gg,\Gd>\to\Z_2$. The two cosets modulo
$\Ker\tilde\kappa$ are represented by~$1$ and~$\Gd$, and the
standard calculation shows that $\Ker\tilde\kappa$ is the free
group generated by $\Ga,\bGa,\Gb,\bGb,\Gg,\bGg,\Gd^2$. The
kernel~$N$ of the epimorphism $\Ker\tilde\kappa\onto\pi_1$ is
normally generated \emph{in $\<\Ga,\Gb,\Gg,\Gd>$} by~$\Gd^2$ and
$R'_j$, $j=1,\ldots,s$. Hence, one can remove
the generator~\smash{$\Gd^2$} from the
presentation. Besides,
since the conjugation by~$\Gd$ is not an
inner automorphism of $\Ker\tilde\kappa$, one should add the
conjugates $\bar R'_j=\Gd^{-1}R'_j\Gd$ to obtain a set normally
generating~$N$ in $\Ker\tilde\kappa$. The resulting presentation
of~$\pi_1$ is the one stated in the lemma.
\endproof

\Remark
Note that $\spbar\:w\mapsto\bar w=\Gd w\Gd$ is an involutive
automorphism of~$\pi_1$. Hence, whenever a relation $R=1$ holds
in~$\pi_1$, the relation $\bar R=1$ also holds.
\endRemark

\subsection{The set of singularities $(3\bE_6)\splus\bA_1$}
\label{s.3e6+a1}
Take for~$\L$ the section given by~\eqref{eq.3e6+a1}. The pair
$(\B,\L)$ looks as shown in Figure~\ref{fig.2e6+2a2+a3}, and the
singular fibers are listed in~\ref{s.cusp}.
The generators
$\zeta_1=\Ga$, $\zeta_2=\Gd$, $\zeta_3=\Gb$, $\zeta_4=\Gg$
for~$\bpi$ are subject
to the relations
$$
\alignat2
&(\Gd\Gb)^2=(\Gb\Gd)^2&\qquad&
 \text{(the tangency point $x=2$)},\\\allowdisplaybreak
&(\Gd\Gb)\Gb(\Gd\Gb)\1=\Gg&&
 \text{(the vertical tangent $x=4$)},\\\allowdisplaybreak
&[\Gd,\Ga\Gd\Gb\Ga]=1,\quad
\Ga\Gd\Gb\Ga=\Gb\Ga\Gd\Gb&&
 \text{(the cusp $x=0$)},\\\allowdisplaybreak
&[\Gd,(\Ga\Gd\Gb)\Gg(\Ga\Gd\Gb)\1]=1&&
 \text{(the transversal intersection $x=-16$)},\\\allowdisplaybreak
&(\Ga\Gd\Gb\Gg)^2=1&&\text{(the relation at infinity)}.
\endalignat
$$
Letting $\Gd^2=1$ and passing to $\Ga$, $\bGa$, $\Gb$, $\bGb$,
$\Gg$, $\bGg$, see Lemma~\ref{bpi->pi}, one can rewrite these
relations in the following form:
$$
\gather
[\Gb,\bGb]=1,
\eqtag\label{eq.1.1}\\\allowdisplaybreak
\Gg=\bGb,\quad \bGg=\Gb,
\eqtag\label{eq.1.2}\\\allowdisplaybreak
\Ga\bGb\bGa=\bGa\Gb\Ga=\Gb\Ga\bGb=\bGb\bGa\Gb,
\eqtag\label{eq.1.3}\\\allowdisplaybreak
\Ga\Gb\Ga\1=\bGa\bGb\bGa\1,
\eqtag\label{eq.1.4}\\\allowdisplaybreak
\Ga\bGb\Gb\bGa\Gb\bGb=1.
\eqtag\label{eq.1.5}
\endgather
$$
(In~\eqref{eq.1.4} and~\eqref{eq.1.5}, we eliminate~$\Gg$
using~\eqref{eq.1.2}.) Now, one can use the last
relation in~\eqref{eq.1.3} to eliminate~$\bGa$: one has
$\bGa=\bGb\1\Gb\Ga\bGb\Gb\1$. Substituting this expression to
$\Ga\bGb\bGa=\Gb\Ga\bGb$ and $\bGa\Gb\Ga=\Gb\Ga\bGb$
in~\eqref{eq.1.3} and
using~\eqref{eq.1.1}, one obtains, respectively, the braid
relations $\Ga\Gb\Ga=\Gb\Ga\Gb$ and $\Ga\bGb\Ga=\bGb\Ga\bGb$.
Conjugating by~$\Gd$, one also has $\bGa\Gb\bGa=\Gb\bGa\Gb$
and $\bGa\bGb\bGa=\bGb\bGa\bGb$. Then, \eqref{eq.1.4} turns to
$\Gb\1\Ga\Gb=\bGb\1\bGa\bGb$ and, eliminating~$\bGa$, one obtains
$[\Ga,\bGb^2\Gb\2]=1$. Finally, eliminating~$\bGa$ from the last
relation~\eqref{eq.1.5}, one gets $\Ga\Gb^2\Ga\bGb^2=1$. Thus, the
map $\Gb\mapsto\Gs_1$, $\Ga\mapsto\Gs_2$, $\bGb\mapsto\Gs_3$
establishes an isomorphism
$$
\pi_1(\Cp2\sminus B)=\BG4/
 \<[\Gs_2,\Gs_1^2\Gs_3\2],\ \Gs_2\Gs_1^2\Gs_2\Gs_3^2>.
$$
It remains to notice that, in the presence of the second relation
in the presentation above, the first one turns into
$[\Gs_2,\Gs_1^2\Gs_2\Gs_1^2\Gs_2]=1$, or
$[\Gs_2,(\Gs_1\Gs_2)^3]=1$, which holds automatically.
Thus, one has
$$
\pi_1(\Cp2\sminus B)=\BG4/\Gs_2\Gs_1^2\Gs_2\Gs_3^2.
\eqtag\label{eq.G0}
$$

\corollary\label{e6->>pi}
Let~$\MB$ be
a Milnor ball about a type~$\bE_6$ singular point of~$B$.
Then the inclusion homomorphism
$\pi_1(\MB\sminus B)\to\pi_1(\Cp2\sminus B)$ is onto.
\endcorollary

\proof
Since any pair of type~$\bE_6$ singular points can be permuted by
a stable symmetry of~$B$, see~\cite{symmetric},
it suffices to prove the
statement for the type~$\bE_6$ point
resulting from the cusp of~$\B$. In
this case, the statement follows from~\eqref{eq.1.2}, as $\Ga$,
$\bGa$, $\Gb$, and~$\bGb$ are all
in the image of $\pi_1(\MB\sminus B)$.
\endproof

\subsection{The set of singularities $(2\bE_6\splus2\bA_2)\splus\bA_3$}
\label{s.2e6+2a2+a3}
Take for~$\L$ the section given by~\eqref{eq.2e6+2a2+a3}. The pair
$(\B,\L)$ is plotted in Figure~\ref{fig.2e6+2a2+a3}, and the
singular fibers are listed in~\ref{s.tangent}.
The generators
$\zeta_1=\Gd$, $\zeta_2=\Ga$, $\zeta_3=\Gb$, $\zeta_4=\Gg$
for~$\bpi$ are subject
to the relations
$$
\alignat2
&[\Gd,\Ga]=1&\qquad&
 \text{(the transversal intersection $x\approx.0459$)},\\\allowdisplaybreak
&\Ga\Gb\Ga=\Gb\Ga\Gb&&
 \text{(the cusp $x=0$)},\\\allowdisplaybreak
&[\Gd,\Gb\Ga\1\Gg\Ga\Gb\1]=1&&
 \text{(the transversal intersection $x\approx-19.1$)},\\\allowdisplaybreak
&(\Gd\Gb)^4=(\Gb\Gd)^4&&
 \text{(the tangency point $x=.625$)},\\\allowdisplaybreak
&(\Gd\Gb)^2\Gb(\Gd\Gb)\2=\Gg&&
 \text{(the vertical tangent $x=4$)},\\\allowdisplaybreak
&(\Gd\Ga\Gb\Gg)^2=1&&\text{(the relation at infinity)}.
\endalignat
$$
(The third relation is simplified using $[\Gd,\Ga]=1$.)
Letting $\Gd^2=1$ and passing to $\Ga=\bGa$,
$\Gb$, $\bGb$, $\Gg$, $\bGg$, see Lemma~\ref{bpi->pi},
one can rewrite these
relations as follows:
$$
\gather
\Ga=\bGa,
\eqtag\label{eq.2.1}\\\allowdisplaybreak
\Ga\Gb\Ga=\Gb\Ga\Gb,\quad \Ga\bGb\Ga=\bGb\Ga\bGb,
\eqtag\label{eq.2.2}\\\allowdisplaybreak
\Gb\Ga\1\bGb\Gb\bGb\1\Ga\Gb\1=\bGb\Ga\1\Gb\bGb\Gb\1\Ga\bGb\1,
\eqtag\label{eq.2.3}\\\allowdisplaybreak
(\bGb\Gb)^2=(\Gb\bGb)^2,
\eqtag\label{eq.2.4}\\\allowdisplaybreak
\bGb\Gb\bGb\1=\Gg,\quad\Gb\bGb\Gb\1=\bGg,
\eqtag\label{eq.2.5}\\\allowdisplaybreak
\Ga\bGb\Gb\bGb\Gb\1\Ga\Gb\bGb\Gb\bGb\1=1.
\eqtag\label{eq.2.6}
\endgather
$$
(We use~\eqref{eq.2.1} and~\eqref{eq.2.5} to eliminate~$\bGa$,
$\Gg$, and~$\bGg$ in the other relations.) Thus,
$$
\pi_1(\Cp2\sminus B)=G_3:=\bigl<\Ga,\Gb,\bGb\bigm|
 \text{\eqref{eq.2.2}--\eqref{eq.2.4}, \eqref{eq.2.6}}\bigr>.
\eqtag\label{eq.G3}
$$

\midinsert
\centerline{\picture{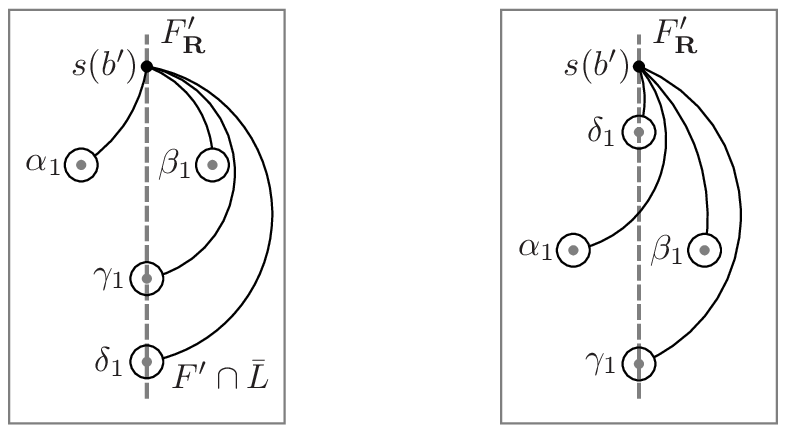}}
\figure
Generators in $F'=\{x=b'=\const\ll0\}$
\endfigure\label{fig.basis2}
\endinsert

The following statement is a consequence of the
monodromy computation.

\lemma\label{e6.a3}
Let $\Ga_1$, $\Gb_1$, $\Gg_1$, $\Gd_1$ be the basis in a fiber
$F'=\{x=\const\ll0\}$ shown in Figure~\ref{fig.basis2}, left.
Then,
considering $\Ga_1$, $\Gb_1$, and $\Gg_1$ as elements of~$\bpi$,
one has
$\Ga_1=\bGb$, $\Gb_1=\Gb\1\Ga\Gb$, and $\Gg_1=\Gg$.
\qed
\endlemma

\corollary\label{e6->>pi.a3}
Let~$\MB$ be
a Milnor ball about a type~$\bE_6$ singular point of~$B$.
Then the inclusion homomorphism
$\pi_1(\MB\sminus B)\to\pi_1(\Cp2\sminus B)$ is onto.
\endcorollary

\proof
In view of~\eqref{eq.2.5}, one has $\Gb=\Ga_1\1\Gg_1\Ga_1$.
Then $\Ga=\Gb\Gb_1\Gb\1$; hence, the elements
$\Ga_1$, $\Gb_1$, and~$\Gg_1$ generate the group. On the other
hand, $\Ga_1$, $\Gb_1$, $\Gg_1$ are in the image of
$\pi_1(\MB\sminus B)$.
\endproof

\subsection{The set of singularities
$(2\bE_6\splus\bA_5)\splus\bA_2$:
the first family}
\label{s.2e6+a5+a2.1}
Take for $\L$ the section given by~\eqref{eq.2e6+a5+a2.1}.
The pair $(\B,\L)$ looks as shown in
Figure~\ref{fig.2e6+2a2+a3}, and the singular fibers are listed
in~\ref{s.inflection}.
The generators
$\zeta_1=\Ga$, $\zeta_2=\Gb$, $\zeta_3=\Gd$, $\zeta_4=\Gg$
satisfy the following
relations:
$$
\alignat2
&[\Gd,\Ga\Gb]=1,\quad \Gd\Ga\Gb\Ga=\Gb\Ga\Gb\Gd
 &&\text{(the cusp $x=0$)},\\\allowdisplaybreak
&(\Gb\Gd)^3=(\Gd\Gb)^3&&
 \text{(the tangency point $x\approx.951$)},\\\allowdisplaybreak
&(\Gb\Gd)\Gb(\Gb\Gd)\1=\Gg&&
 \text{(the vertical tangent $x=4$)},\\\allowdisplaybreak
&[\Gd,\Ga\1\Gg\Ga]=1&&
 \text{(the transversal intersection $x\approx-18.4$)},\\\allowdisplaybreak
&(\Ga\Gb\Gd\Gg)^2=1&\qquad&\text{(the relation at infinity)}.
\endalignat
$$
Letting $\Gd^2=1$ and passing to $\Ga$,
$\bGa$, $\Gb$, $\bGb$, $\Gg$, $\bGg$, see Lemma~\ref{bpi->pi},
one obtains
$$
\gather
\Ga\Gb=\bGa\bGb,\quad\bGa\bGb\bGa=\Gb\Ga\Gb,\quad
 \Ga\Gb\Ga=\bGb\bGa\bGb,
\eqtag\label{eq.3.1}\\\allowdisplaybreak
\bGb\Gb\bGb=\Gb\bGb\Gb,
\eqtag\label{eq.3.2}\\\allowdisplaybreak
\Gb\bGb\Gb\1=\Gg,\quad\bGb\Gb\bGb\1=\bGg,
\eqtag\label{eq.3.3}\\\allowdisplaybreak
\Ga\1\Gg\Ga=\bGa\1\bGg\bGa,
\eqtag\label{eq.3.4}\\\allowdisplaybreak
\Ga\Gb\bGg\Ga\Gb\Gg=1.
\eqtag\label{eq.3.5}
\endgather
$$
The cusp relations~\eqref{eq.3.1} can be rewritten in the form
$$
\bGa=(\Ga\Gb)\1\Gb(\Ga\Gb),\quad
\bGb=(\Ga\Gb)\Ga(\Ga\Gb)\1,\quad
(\Ga\Gb)^3=(\Gb\Ga)^3,
\eqtag\label{eq.3.6}
$$
or, in terms of~$\bGa$, $\bGb$, in the form
$$
\Ga=(\bGa\bGb)\1\bGb(\bGa\bGb),\quad
\Gb=(\bGa\bGb)\bGa(\bGa\bGb)\1,\quad
(\bGa\bGb)^3=(\bGb\bGa)^3.
\eqtag\label{eq.3.7}
$$
Geometrically, one has
$\pi_1(\MB\sminus B)=\bigl<\Ga,\Gb\bigm|(\Ga\Gb)^3=(\Gb\Ga)^3\bigr>$,
where $\MB$ is a Milnor ball around the type~$\bA_5$ singular
point.

Writing~\eqref{eq.3.5} as $\Ga\Gb\bGg\bGa\bGb\Gg=1$ and
eliminating~$\Gg$ and~$\bGg$ using~\eqref{eq.3.3}
and~\eqref{eq.3.2}, we can rewrite this relation in the form
$$
\Ga\bGb\Gb\bGa\Gb\bGb=1.
\eqtag\label{eq.3.8}
$$
Eliminating~$\Gg$ and~$\bGg$ from~\eqref{eq.3.4}, we obtain
$$
\Ga\1\Gb\bGb\Gb\1\Ga=\bGa\1\bGb\Gb\bGb\1\bGa.
\eqtag\label{eq.3.9}
$$
Thus, we have
$$
\pi_1(\Cp2\sminus B)=G_2':=\bigl<\Ga,\Gb\bigm|
 \text{$(\Ga\Gb)^3=(\Gb\Ga)^3$, \eqref{eq.3.2}, \eqref{eq.3.8},
 \eqref{eq.3.9}}\bigr>,
\eqtag\label{eq.G2.1}
$$
where $\bGa$ and $\bGb$ are the words given by~\eqref{eq.3.6}. I
could not find any substantial simplification of this
presentation. An alternative presentation of~$G_2'$ (as well as
of the group~$G_2''$ introduced in~\eqref{eq.G2.2} below)
is given in Eyral, Oka~\cite{EyralOka}.

As a part of computing the braid monodromy, we get the following
lemma.

\lemma\label{e6.a5.1}
Let $\Ga_1$, $\Gb_1$, $\Gg_1$, $\Gd_1$ be the basis in a fiber
$F'=\{x=\const\ll0\}$ shown in Figure~\ref{fig.basis2}, left.
Then,
considering $\Ga_1$, $\Gb_1$, and $\Gg_1$ as elements of~$\bpi$,
one has
$\Ga_1=\Gb$, $\Gb_1=\bGb\1\bGa\bGb$, and $\Gg_1=\Gg$.
\qed
\endlemma

\corollary\label{e6->>pi.a5.1}
Let~$\MB$ be
a Milnor ball about a type~$\bE_6$ singular point of~$B$.
Then the inclusion homomorphism
$\pi_1(\MB\sminus B)\to\pi_1(\Cp2\sminus B)$ is onto.
\endcorollary

\proof
Due to~\eqref{eq.3.3}, one has $\bGb=\Ga_1\1\Gg_1\Ga_1$. Then
$\bGa=\bGb\Gb_1\bGb\1$ and, in view of~\eqref{eq.3.7}
and~\eqref{eq.3.3}, $\bGa$ and $\bGb$ generate the group.
\endproof

\subsection{The set of singularities
$(2\bE_6\splus\bA_5)\splus\bA_2$:
the second family}
\label{s.2e6+a5+a2.2}
Now, let~$\L$ be the section given by~\eqref{eq.2e6+a5+a2.2}. The
pair $(\B,\L)$ is plotted in Figure~\ref{fig.2e6+a5+a2}, and the
singular fibers are listed in~\ref{s.inflection}.
The generators for~$\pi_F$ are
$\zeta_1=\Ga$, $\zeta_2=\Gb$, $\zeta_3=\Gd$, $\zeta_4=\Gg$,
and the
relations are:
$$
\alignat2
&[\Gd,\Ga\Gb]=1,\quad \Gd\Ga\Gb\Ga=\Gb\Ga\Gb\Gd
 &&\text{(the cusp $x=0$)},\\\allowdisplaybreak
&(\Gg\Gd)^3=(\Gd\Gg)^3&&
 \text{(the tangency point $x\approx3.55$)},\\\allowdisplaybreak
&(\Gd\Gg\Gd)\Gg(\Gd\Gg\Gd)\1=\Gb&&
 \text{(the vertical tangent $x=4$)},\\\allowdisplaybreak
&[\Gd,\Gg\Ga\Gg\1]=1&&
 \text{(the transversal intersection $x\approx4.94$)},\\\allowdisplaybreak
&(\Ga\Gb\Gd\Gg)^2=1&\qquad&\text{(the relation at infinity)}.
\endalignat
$$
Let $\Gd^2=1$ and pass to the generators $\Ga$,
$\bGa$, $\Gb$, $\bGb$, $\Gg$, $\bGg$, see Lemma~\ref{bpi->pi}.
Then,
in addition to the cusp relations~\eqref{eq.3.6}
(or~\eqref{eq.3.1}\,) and relation at infinity~\eqref{eq.3.5},
we obtain
$$
\gather
\Gg\bGg\Gg=\bGg\Gg\bGg,
\eqtag\label{eq.4.1}\\\allowdisplaybreak
\bGg\Gg\bGg\1=\Gb,\quad\Gg\bGg\Gg\1=\bGb
\eqtag\label{eq.4.2}\\\allowdisplaybreak
\Gg\Ga\Gg\1=\bGg\bGa\bGg\1.
\eqtag\label{eq.4.3}
\endgather
$$
Thus,
$$
\pi_1(\Cp2\sminus B)=G_2'':=\bigl<\Ga,\Gb,\Gg,\bGg\bigm|
 \text{$(\Ga\Gb)^3=(\Gb\Ga)^3$, \eqref{eq.3.5},
 \eqref{eq.4.1}--\eqref{eq.4.3}}\bigr>,
\eqtag\label{eq.G2.2}
$$
where $\bGa$ and $\bGb$ are the words given by~\eqref{eq.3.6}.
Note that one can eliminate either~$\bGg$, using~\eqref{eq.3.5},
or~$\Gb$, using~\eqref{eq.4.2}.

Extending the braid monodromy beyond the cusp of~$B$ (to the
negative values of~$x$), we obtain the
following statement.

\lemma\label{e6.a5.2}
Let $\Gd_1$, $\Ga_1$, $\Gb_1$, $\Gg_1$ be the basis in a fiber
$F'=\{x=\const\ll0\}$ shown in Figure~\ref{fig.basis2}, right.
Then,
considering $\Ga_1$, $\Gb_1$, and $\Gg_1$ as elements of~$\bpi$,
one has
$\Ga_1=\bGb$, $\Gb_1=\bGb\1\bGa\bGb$, and $\Gg_1=\Gg$.
\qed
\endlemma

\corollary\label{e6->>pi.a5.2}
Let~$\MB$ be
a Milnor ball about a type~$\bE_6$ singular point of~$B$.
Then the inclusion homomorphism
$\pi_1(\MB\sminus B)\to\pi_1(\Cp2\sminus B)$ is onto.
\endcorollary

\proof
In view of~\eqref{eq.3.7} and~\eqref{eq.3.5}, the elements
$\bGa=\Ga_1\Gb_1\Ga_1\1$, $\bGb=\Ga_1$,
and $\Gg=\Gg_1$ generate the group.
\endproof

\subsection{Comparing the two groups}
Let~$B'$ and~$B''$ be the sextics considered
in~\ref{s.2e6+a5+a2.1} and~\ref{s.2e6+a5+a2.2}, respectively, so
that their fundamental groups are~$G_2'$ and~$G_2''$. As explained
in Eyral, Oka~\cite{EyralOka}, the profinite completions of~$G_2'$
and~$G_2''$ are isomorphic (as the two curves are conjugate over
an algebraic number field). Whether $G_2'$ and~$G_2''$ themselves
are isomorphic is still an open question. Below, we suggest an
attempt to distinguish the two groups geometrically.

\proposition\label{a5->>pi}
Let~$\MB$ be
a Milnor ball about the type~$\bA_5$ singular point of~$B'$.
Then the inclusion homomorphism
$\pi_1(\MB\sminus B')\to\pi_1(\Cp2\sminus B')$ is onto.
\endproposition

\proof
According to~\eqref{eq.G2.1}, the group
$\pi_1(\Cp2\sminus B')=G_2'$ is generated by $\Ga$ and~$\Gb$,
which are both in the image of $\pi_1(\MB\sminus B')$.
\endproof

\conjecture\label{a5not->>pi}
Let~$\MB$ be
a Milnor ball about the type~$\bA_5$ singular point of~$B''$.
Then the image of the inclusion homomorphism
$\pi_1(\MB\sminus B'')\to\pi_1(\Cp2\sminus B'')$ does not
contain~$\Gg$ or~$\bGg$.
\endconjecture

\Remark
If true, Conjecture~\ref{a5not->>pi} together with
Proposition~\ref{a5->>pi} would provide a topological distinction
between pairs $(\Cp2,B')$ and $(\Cp2,B'')$. Note that,
according to~\cite{JAG}, the two pairs are not diffeomorphic.
\endRemark

\subsection{Other symmetric sets of singularities}
\label{s.stable}
The set of singularities $(3\bE_6)$ is obtained by perturbing~$\L$
in Section~\ref{s.3e6+a1} to a section tangent to~$\B$ at the cusp
and transversal to~$\B$ otherwise. This procedure
replaces~\eqref{eq.1.1} with $\bGb=\Gb$ or, alternatively,
introduces a relation $\Gs_3=\Gs_1$ in~\eqref{eq.G0}. The
resulting group is $\BG3/(\Gs_1\Gs_2)^3$.

The sets of singularities of the form
$(2\bE_6\splus2\bA_2)\splus{\ldots}$ are obtained by
perturbing $\L$ in Section~\ref{s.2e6+2a2+a3}.
If $\L$ is perturbed to a double tangent (the set of singularities
$(2\bE_6\splus2\bA_2)\splus2\bA_1$), relation~\eqref{eq.2.4}
is replaced with $[\Gb,\bGb]=1$. Then, \eqref{eq.2.6} turns to
$\Ga\bGb^2\Ga\Gb^2=1$, and \eqref{eq.2.3} turns to
$$
\Gb\underline{\Ga\1\Gb\Ga}\Gb\1=\bGb\underline{\Ga\1\bGb\Ga}\bGb\1.
$$
Replacing the underlined expressions using the braid
relations~\eqref{eq.2.2} converts this relation to
$\Gb^2\Ga\Gb\2=\bGb^2\Ga\bGb\2$, \ie, $[\Ga,\bGb^2\Gb\2]=1$.
As explained in~\ref{s.3e6+a1},
the map
$\Gb\mapsto\Gs_1$, $\Ga\mapsto\Gs_2$, $\bGb\mapsto\Gs_3$
establishes an isomorphism
$\pi_1(\Cp2\sminus B)=\BG4/\Gs_2\Gs_1^2\Gs_2\Gs_3^2$.

Any other perturbation of~$\L$ produces an extra point of
transversal intersection with~$\B$, replacing~\eqref{eq.2.4} with
$\Gb=\bGb$. The resulting group is $\BG3/(\Gs_1\Gs_2)^3$.

Finally, the sets of singularities $(2\bE_6\splus\bA_5)\splus\bA_1$
and $(2\bE_6\splus\bA_5)$ are obtained by perturbing the inflection
tangency point of~$\L$ and~$\B$ in Section~\ref{s.2e6+a5+a2.1}.
This procedure replaces~\eqref{eq.3.2} with $\bGb=\Gb$. Then, from
the first relation in~\eqref{eq.3.1} one has $\bGa=\Ga$,
relation~\eqref{eq.3.3} results in $\Gg=\bGb=\Gb$,
and relation~\eqref{eq.3.5} turns to $(\Ga\Gb^2)^2=1$.
Hence, the group is $\BG3/(\Gs_1\Gs_2)^3$. (Note that
$(\Gs_1\Gs_2^2)^2=(\Gs_1\Gs_2)^3$ in $\BG3$.)

\subsection{Proof of Theorem~\ref{th.proper}}\label{proof.proper}
The fact that the perturbation epimorphisms
$G_2',G_2''\onto\BG3/(\Gs_1\Gs_2)^3$ are proper is proved in Eyral,
Oka~\cite{EyralOka}, where it is shown that the Alexander module
of a sextic with the set of singularities
$(2\bE_6\splus\bA_5)\splus\bA_2$ has a
torsion summand $\Z_2\times\Z_2$, whereas the Alexander modules
of all other groups
listed in Theorem~\ref{th.group}
can easily be shown to be $\Z[t]/(t^2-t+1)$. (In other words, the
abelianization of the commutant of~$G_2'$ or $G_2''$ is equal to
$\Z_2\times\Z_2\times\Z\times\Z$, and for all other groups it
equals $\Z\times\Z$.)

The epimorphism
$$
\Gf_0\:G_0=\BG4/\Gs_2\Gs_1^2\Gs_2\Gs_3^2\onto\BG3/(\Gs_1\Gs_2)^3
$$
is considered in Oka, Pho~\cite{OkaPho}. One can observe that both
braids $\Gs_2\Gs_1^2\Gs_2\Gs_3^2$ and $(\Gs_1\Gs_2)^3$ in the
definition of the groups are pure, \ie, belong to the kernels of
the respective canonical epimorphism
$\BG{n}\onto\BG{n}/\Gs_1^2=\SG{n}$. Furthermore, $\Gf_0$ takes each
of the standard generators $\Gs_1$, $\Gs_2$, $\Gs_3$ of~$\BG4$ to
a conjugate of $\Gs_1$. Hence, the induced epimorphism
$G_0/\Gf^{-1}(\Gs_1^2)=\SG4\onto\BG3/\Gs_1^2=\SG3$ is proper, and
so is~$\Gf_0$.

A similar argument applies to the epimorphism
$\Gf_3\:G_3\onto G_0$, which takes each generator $\Ga$, $\Gb$,
$\bGb$ of~$G_3$ to a conjugate of~$\Gs_1\in G_0$. The induced
epimorphism
$$
G_3/\Gf_3^{-1}(\Gs_1^2)=\SL(2,\Bbb F_3)\onto
 G_0/\Gs_1^2=\SG4=\PSL(2,\Bbb F_3)
$$
is proper; hence, so is $\Gf_3$. (Alternatively, one can compare
$G_3/\Gf_3^{-1}(\Gs_1^4)$ and $G_0/\Gs_1^4$, which are finite
groups of order $3\cdot2^9$ and $3\cdot2^6$, respectively. The finite
quotients of~$G_3$ and~$G_0$ were computed using {\tt GAP}~\cite{GAP}.)
\qed

\section{Perturbations\label{S.perturbations}}

\subsection{Perturbing a singular point}\label{s.perturbations}
Consider a singular point~$P$ of a plain curve $B$ and a Milnor
ball~$\MB$ around~$P$. Let~$B'$ be a nontrivial (\ie, not
equisingular) perturbation of~$B$ such that,
during the perturbation, the curve remains transversal to
$\partial\MB$.

\lemma\label{pert.E6}
In the notation above, let $P$ be of type~$\bE_6$. Then $B'\cap\MB$
has one of the following sets of singularities\rom:
\roster
\item\local{E6.max}
$2\bA_2\splus\bA_1$\rom: one has $\pi_1(\MB\sminus B')=\BG4$\rom;
\item\local{E6.braid}
$\bA_5$ or $2\bA_2$\rom: one has $\pi_1(\MB\sminus B')=\BG3$\rom;
\item\local{E6.cyclic}
$\bD_5$, $\bD_4$, $\bA_4\splus\bA_1$, $\bA_4$, $\bA_3\splus\bA_1$,
$\bA_3$, $\bA_2\splus k\bA_1$ \rom($k=0$, $1$, or~$2$\rom), or
$k\bA_1$ \rom($k=0$, $1$, $2$, or~$3$\rom)\rom:
one has $\pi_1(\MB\sminus B')=\Z$.
\endroster
\endlemma

\proof
The perturbations of a simple singularity are enumerated by the
subgraphs of its Dynkin graph, see E.~Brieskorn~\cite{Brieskorn}
or G.~Tjurina~\cite{Tjurina}.
For the fundamental group, observe that
the space $\MB\sminus B$ is diffeomorphic to
$\Cp2\sminus(C\cup L)$, where $C\subset\Cp2$ is a plane quartic
with a type~$\bE_6$ singular point, and $L$ is a line with a
single quadruple intersection point with~$C$. Then, the
perturbations of~$B$ inside~$\MB$ can be regarded as perturbations
of~$C$ keeping the point of quadruple intersection with~$L$,
see~\cite{quintics}, and the perturbed fundamental group
$\pi_1(\Cp2\sminus(C'\cup L)\cong\pi_1(\MB\sminus B')$
is found in~\cite{groups}.
\endproof

\lemma\label{pert.A5}
In the notation above, let $P$ be of type~$\bA_5$. Then $B'\cap\MB$
has one of the following sets of singularities\rom:
\roster
\item\local{A5.braid}
$2\bA_2$\rom: one has $\pi_1(\MB\sminus B')=\BG3$\rom;
\item\local{A5.abelian}
$\bA_3\splus\bA_1$ or $3\bA_1$\rom:
one has $\pi_1(\MB\sminus B')=\Z\times\Z$\rom;
\item\local{A5.cyclic}
$\bA_4$, $\bA_3$, $\bA_2\splus\bA_1$, $\bA_2$, or
$k\bA_1$ \rom($k=0$, $1$, or~$2$\rom)\rom:
one has $\pi_1(\MB\sminus B')=\Z$.
\endroster
\endlemma

\lemma\label{pert.A2}
In the notation above, let $P$ be of type~$\bA_2$. Then $B'\cap\MB$
has the set of singularities $\bA_1$ or $\varnothing$, and
one has $\pi_1(\MB\sminus B')=\Z$.
\endlemma

\proof[Proof of Lemmas~\ref{pert.A5} and~\ref{pert.A2}]
Both statements are a well known property of type~$\bA$ singular
points: any perturbation of a type~$\bA_p$ singular point has the
set of singularities $\bigsplus\bA_{p_i}$ with
$d=(p+1)-\sum(p_i+1)\ge0$, and the group $\pi_1(\MB\sminus B')$ is
given by $\<\Ga,\Gb\,|\,\Gs^s\Ga=\Ga,\ \Gs^s\Gb=\Gb>$, where $\Gs$
is the standard generator of the braid group~$\BG2$ acting on
$\<\Ga,\Gb>$ and $s=1$ if $d>0$ or $s=\gcd(p_i+1)$ if $d=0$.
\endproof

\proposition\label{E6.onto}
Let~$B$ be a plane sextic of torus type with at least two
type~$\bE_6$ singularities, and let $\MB$ be a Milnor ball about a
type~$\bE_6$ singular point of~$B$. Then the inclusion
homomorphism $\pi_1(\MB\sminus B)\to\pi_1(\Cp2\sminus B)$
is onto.
\endproposition

\proof
The proposition is an immediate consequence of
Corollaries~\ref{e6->>pi}, \ref{e6->>pi.a3},
\ref{e6->>pi.a5.1} and~\ref{e6->>pi.a5.2}.
\endproof

\corollary\label{E6.perturbed}
Let~$B$ be a plane sextic of torus type with at least two
type~$\bE_6$ singular points, and let $B'$ be a perturbation of~$B$.
\roster
\item\local{E6.perturbed.abelian}
If at least one of the type~$\bE_6$ singular points of~$B$ is
perturbed as in~\iref{pert.E6}{E6.cyclic}, then
$\pi_1(\Cp2\sminus B')=\CG6$.
\item\local{E6.perturbed.braid}
If at least one of the type~$\bE_6$ singular points of~$B$ is
perturbed as in~\iref{pert.E6}{E6.braid}
and $B'$ is still of torus type,
then $\pi_1(\Cp2\sminus B')=\BG3/(\Gs_1\Gs_2)^3$.
\endroster
\endcorollary

\proof
Let~$\MB$ be a Milnor ball about the type~$\bE_6$ singular point in
question. Due to
Proposition~\ref{E6.onto},
the inclusion homomorphism
$\pi_1(\MB\sminus B)\to\pi_1(\Cp2\sminus B)$ is onto. Hence, in
case~\loccit{E6.perturbed.abelian}, there
is an epimorphism $\Z\onto\pi_1(\Cp2\sminus B')$, and
in case~\loccit{E6.perturbed.braid}, there
is an epimorphism $\BG3\onto\pi_1(\Cp2\sminus B')$.
In the former case, the epimorphism above implies that the group
is abelian, hence~$\CG6$.
In the latter case, the central element $(\Gs_1\Gs_2)^3\in\BG3$
projects to $6\in\Z=\BG3/[\BG3,\BG3]$; since the abelianization of
$\pi_1(\Cp2\sminus B')$ is $\CG6$, the epimorphism above must
factor through an epimorphism
$G:=\BG3/(\Gs_1\Gs_2)^3\onto\pi_1(\Cp2\sminus B')$. On the other
hand, since $B'$ is assumed to be of torus type, there is an
epimorphism $\pi_1(\Cp2\sminus B')\onto G$, and as
$G\cong\PSL(2,\ZZ)$ is Hopfian (as it is obviously residually finite),
each of the two epimorphisms is bijective.
\endproof

\corollary\label{A5.perturbed}
Let~$B$ be a plane sextic as in~\ref{s.2e6+a5+a2.1},
and let $B'$ be a perturbation of~$B$ such that
the type~$\bA_5$ singular point is
perturbed as in~\iref{pert.A5}{A5.abelian} or~\ditto{A5.cyclic}.
Then one has $\pi_1(\Cp2\sminus B')=\CG6$.
\endcorollary

\proof
Due to Proposition~\ref{a5->>pi} and
Lemma~\ref{pert.A5}, the group of the perturbed sextic~$B'$ is
abelian. Since $B'$ is irreducible, $\pi_1(\Cp2\sminus B')=\CG6$.
\endproof

\corollary\label{A2.perturbed}
Let~$B$ be a plane sextic as in~\ref{s.2e6+2a2+a3},
and let $B'$ be a perturbation of~$B$ such that
an inner type~$\bA_2$ singular point of~$B$ is
perturbed to~$\bA_1$ or~$\varnothing$.
Then one has $\pi_1(\Cp2\sminus B')=\CG6$.
\endcorollary

\proof
Let~$P$ be the inner type~$\bA_2$ singular point perturbed, and
let $\MB$ be a Milnor ball about~$P$. In the notation of
Section~\ref{s.2e6+2a2+a3}, the group $\pi_1(\MB\sminus B)$ is
generated by~$\Ga$ and~$\Gb$ (or $\bGa=\Ga$ and~$\bGb$ for the
other point), and the perturbation results in an extra relation
$\Ga=\Gb$. Then \eqref{eq.2.3} implies $\bGb=\Gb$ and the group is
cyclic.
\endproof

%
%
%

\midinsert
\table\label{tab.nontorus}
Sextics with abelian fundamental group
\endtable
\hbox to\hsize{\hss
\TAB\NonTorus
2 e[6] + d[5] + a[1],               [37, 4], [2, 2, 0], 18, [3, 2, false]
2 e[6] + a[4] + 2 a[1],             [38, 4], [2, 2, 1], 18, [3, 2, false]
2 e[6] + a[3] + a[2] + a[1],        [40, 5], [2, 3, 0], 18, [3, 3, false]
2 e[6] + 2 a[2] + a[1],             [38, 6], [1, 4, 0], 17, [4, 4, false]
2 e[6] + a[2] + 3 a[1],             [36, 5], [3, 3, 0], 17, [4, 3, false]
e[6] + 2 d[5] + a[1],               [36, 2], [3, 1, 0], 17, [4, 3, false]
e[6] + d[5] + a[5] + a[1],          [37, 4], [3, 2, 0], 17, [4, 3, false]
e[6] + d[5] + a[4] + 2 a[1],        [37, 2], [3, 1, 1], 17, [4, 3, false]
e[6] + d[5] + a[3] + 2 a[2],        [43, 4], [2, 3, 0], 18, [3, 3, false]
e[6] + d[5] + 3 a[2],               [41, 5], [1, 4, 0], 17, [4, 4, false]
e[6] + d[5] + 2 a[2] + 2 a[1],      [39, 4], [3, 3, 0], 17, [4, 3, false]
e[6] + d[4] + a[3] + 2 a[2],        [41, 4], [3, 3, 0], 17, [4, 3, false]
e[6] + d[4] + 3 a[2],               [39, 5], [2, 4, 0], 16, [5, 4, false]
e[6] + d[4] + 2 a[2] + 2 a[1],      [37, 4], [4, 3, 0], 16, [5, 4, false]
e[6] + a[5] + a[4] + a[1] + a[2],   [42, 5], [2, 3, 1], 18, [3, 3, false]
e[6] + a[5] + a[3] + a[2] + a[1],   [40, 5], [3, 3, 0], 17, [4, 3, false]
e[6] + 2 a[4] + 3 a[1],             [38, 2], [3, 1, 2], 17, [4, 3, false]
e[6] + a[4] + a[3] + 2 a[2] + a[1], [44, 4], [2, 3, 1], 18, [3, 3, false]
e[6] + a[4] + 2 a[2] + 3 a[1],      [40, 4], [3, 3, 1], 17, [4, 3, false]
e[6] + 2 a[3] + 2 a[2] + a[1],      [42, 4], [3, 3, 0], 17, [4, 3, false]
2 d[5] + a[5] + a[1],               [36, 2], [4, 1, 0], 16, [5, 4, false]
2 d[5] + a[4] + 2 a[1],             [36, 0], [4, 0, 1], 16, [5, 4, false]
2 d[5] + a[3] + 2 a[2],             [42, 2], [3, 2, 0], 17, [4, 3, false]
2 d[5] + 3 a[2],                    [40, 3], [2, 3, 0], 16, [5, 3, true]
2 d[5] + 2 a[2] + 2 a[1],           [38, 2], [4, 2, 0], 16, [5, 4, false]
d[5] + d[4] + a[3] + 2 a[2],        [40, 2], [4, 2, 0], 16, [5, 4, false]\ENDTAB
\hss
\TAB\NonTorus
d[5] + d[4] + 3 a[2],               [38, 3], [3, 3, 0], 15, [6, 3, true]
d[5] + d[4] + 2 a[2] + 2 a[1],      [36, 2], [5, 2, 0], 15, [6, 5, false]
d[5] + 2 a[5] + a[1],               [37, 4], [4, 2, 0], 16, [5, 4, false]
d[5] + a[5] + a[4] + 2 a[1],        [37, 2], [4, 1, 1], 16, [5, 4, false]
d[5] + a[5] + a[3] + 2 a[2],        [43, 4], [3, 3, 0], 17, [4, 3, false]
d[5] + 2 a[4] + 3 a[1],             [37, 0], [4, 0, 2], 16, [5, 4, false]
d[5] + a[4] + a[3] + 2 a[2] + a[1], [43, 2], [3, 2, 1], 17, [4, 3, false]
d[5] + a[4] + 2 a[2] + 3 a[1],      [39, 2], [4, 2, 1], 16, [5, 4, false]
d[5] + 2 a[3] + 2 a[2] + a[1],      [41, 2], [4, 2, 0], 16, [5, 4, false]
2 d[4] + a[3] + 2 a[2],             [38, 2], [5, 2, 0], 15, [6, 5, false]
2 d[4] + 3 a[2],                    [36, 3], [4, 3, 0], 14, [7, 4, true]
d[4] + a[5] + a[3] + 2 a[2],        [41, 4], [4, 3, 0], 16, [5, 4, false]
d[4] + a[4] + a[3] + 2 a[2] + a[1], [41, 2], [4, 2, 1], 16, [5, 4, false]
d[4] + a[4] + 2 a[2] + 3 a[1],      [37, 2], [5, 2, 1], 15, [6, 5, false]
d[4] + 2 a[3] + 2 a[2] + a[1],      [39, 2], [5, 2, 0], 15, [6, 5, false]
2 a[5] + a[4] + a[1] + a[2],        [42, 5], [3, 3, 1], 17, [4, 3, false]
2 a[5] + a[3] + a[2] + a[1],        [40, 5], [4, 3, 0], 16, [5, 4, false]
a[5] + 2 a[4] + 2 a[1] + a[2],      [42, 3], [3, 2, 2], 17, [4, 3, false]
a[5] + a[4] + a[3] + 2 a[2] + a[1], [44, 4], [3, 3, 1], 17, [4, 3, false]
a[5] + 2 a[3] + 2 a[2] + a[1],      [42, 4], [4, 3, 0], 16, [5, 4, false]
3 a[4] + 4 a[1],                    [38, 0], [4, 0, 3], 16, [5, 4, false]
2 a[4] + a[3] + 2 a[2] + 2 a[1],    [44, 2], [3, 2, 2], 17, [4, 3, false]
2 a[4] + 2 a[2] + 4 a[1],           [40, 2], [4, 2, 2], 16, [5, 4, false]
a[4] + 2 a[3] + 2 a[2] + 2 a[1],    [42, 2], [4, 2, 1], 16, [5, 4, false]
3 a[3]+2 a[2]+2 a[1],               [40, 2], [5, 2, 0], 15, [6, 5, false]\ENDTAB\hss}
\endinsert

\subsection{Abelian perturbations}
Theorem~\ref{th.nontorus} below lists the sets of singularities
obtained by perturbing at least one inner singular point from a
set listed in Table~\ref{tab.list}, not covered by Nori's
theorem~\cite{Nori}, and not appearing in~\cite{degt.8a2}.

\theorem\label{th.nontorus}
Let $\Sigma$ be a set of singularities obtained from one of those
listed in Table~\ref{tab.nontorus} by several \rom(possibly
none\rom) perturbations $\bA_2\to\bA_1,\varnothing$ or
$\bA_1\to\varnothing$.
Then $\Sigma$ is realized by an irreducible plane
sextic, not of torus type, whose fundamental group is~$\CG6$.
\endtheorem

Altogether, perturbations as in Theorem~\ref{th.nontorus} produce
$244$ sets of singularities not covered by Nori's theorem; $117$
of them are new as compared to~\cite{degt.8a2}.

\proof
Each set of singularities in question is obtained by a perturbation
from one of the sets of singularities listed in
Table~\ref{tab.list}. Furthermore, the perturbation can be chosen
so that at least one type~$\bE_6$ singular point is perturbed as
in~\iref{pert.E6}{E6.cyclic}, or the type~$\bA_5$ singular point
is perturbed as in~\iref{pert.A5}{A5.cyclic}, or at least one
inner cusp is perturbed to~$\bA_1$ or~$\varnothing$. According
to~\cite{degt.8a2}, any such (formal) perturbation is realized by
a family of sextics, and due to
Corollaries~\iref{E6.perturbed}{E6.perturbed.abelian},
\ref{A5.perturbed}, and~\ref{A2.perturbed}, the perturbed sextic
has abelian fundamental group.
\endproof

\subsection{Non-abelian perturbations}
In this section, we treat the few perturbations of torus type that
can be obtained from Table~\ref{tab.list} and do not appear
in~\cite{degt.8a2}.

\theorem\label{th.torus}
Each of the eight sets of singularities listed in
Table~\ref{tab.torus} is realized by an irreducible plane sextic
of torus type whose fundamental group is $\BG3/(\Gs_1\Gs_2)^3$.
\midinsert
\table\label{tab.torus}
Sextics of torus type
\endtable
\hbox to\hsize{\hss
\Tab
(\bE_6\splus2\bA_5)\splus\bA_2\cr
(\bE_6\splus2\bA_5)\splus\bA_1\cr
(\bE_6\splus2\bA_5)\cr
(3\bA_5)\splus\bA_2
\ENDTAB\hss\Tab
(3\bA_5)\splus\bA_1\cr
(3\bA_5)\cr
(\bE_6\splus\bA_5\splus2\bA_2)\splus\bA_3\cr
(2\bA_5\splus2\bA_2)\splus\bA_3
\ENDTAB
\hss}
\endinsert
\endtheorem

Theorem~\ref{th.torus} covers two tame sextics:
$(\bE_6\splus2\bA_5)$ and $(3\bA_5)$. The fundamental groups of these
curves were first found in Oka, Pho~\cite{OkaPho}.

\proof
As in the previous section, we perturb one of the sets of
singularities listed in Table~\ref{tab.list}, this time making
sure that
\roster
\item\local{pert.1}
each type~$\bE_6$ singular point is perturbed as
in~\iref{pert.E6}{E6.max} or~\ditto{E6.braid} (or is not
perturbed at all),
\item\local{pert.2}
each type~$\bA_5$ singular point is perturbed as
in~\iref{pert.A5}{A5.braid} (or is not
perturbed at all),
\item\local{pert.3}
none of the inner cusps is perturbed, and
\item\local{pert.4}
at least one type~$\bE_6$ singular point is perturbed as
in~\iref{pert.E6}{E6.braid}.
\endroster
(Note that, in the case under consideration, inner are the cusps
appearing from the cusp of~$\B$.)
From the arithmetic description of curves of torus type given
in~\cite{degt.Oka} (see also~\cite{JAG}) it follows that any
perturbation satisfying~\loccit{pert.1}--\loccit{pert.3} above
preserves the torus structure; then, in view of~\loccit{pert.4},
Corollary~\iref{E6.perturbed}{E6.perturbed.braid} implies that the
resulting fundamental group is $\BG3/(\Gs_1\Gs_2)^3$.
\endproof

\Remark\label{not.new}
If all type~$\bE_6$ singular points are perturbed as
in~\iref{pert.E6}{E6.max} (or not perturbed at all), the study of
the fundamental group would require more work; in particular, one
would need an explicit description of the homomorphism
$\pi_1(\MB_{\bE_6}\sminus B)\onto\pi_1(\MB_{\bE_6}\sminus B')$. On
the other hand, it is easy to show that such perturbations do not
give anything new compared to~\cite{degt.8a2}. (In fact,
using~\cite{JAG}, one can even show that the deformation classes
of the sextics obtained are the same; it suffices to prove the
connectedness of the deformation families realizing the
sets of singularities
$(\bE_6\splus\bA_5\splus2\bA_2)\splus\bA_2\splus\bA_1$ and
$(\bE_6\splus4\bA_2)\splus\bA_3\splus\bA_1$,
which are maximal in the context.) For this reason, we do not
consider these perturbations here.
\endRemark

\widestnumber\key{EO1}
\refstyle{C}


\widestnumber\no{99}
\Refs

\ref{Bri}
\by E.~Brieskorn
\paper Singular elements of semi-simple algebraic groups
\inbook Actes Congr\`es Inter. Math., Nice
\issue 2
\yr 1970
\pages 279-284
\endref\label{Brieskorn}

\ref{D1}
\by A.~Degtyarev
\paper Isotopy classification of complex plane projective curves of
degree~$5$
\jour Algebra i Analis
\yr 1989
\vol    1
\issue  4
\pages  78--101
\lang Russian
\moreref\nofrills English transl. in
\jour Leningrad Math.~J.
\vol 1
\yr 1990
\issue 4
\pages 881--904
\endref\label{quintics}

\ref{D2}
\by A.~Degtyarev
\paper Quintics in $\C\roman{p}^2$ with nonabelian fundamental group
\jour Algebra i Analis
\yr 1999
\vol    11
\issue  5
\pages  130--151
\lang Russian
\moreref\nofrills English transl. in
\jour Leningrad Math.~J.
\vol 11
\yr 2000
\issue 5
\pages 809--826
\endref\label{groups}

\ref{D3}
\by A.~Degtyarev
\paper On deformations of singular plane sextics
\jour J. Algebraic Geom.
\vol 17
\yr 2008
\pages 101--135
\endref\label{JAG}

\ref{D4}
\by A.~Degtyarev
\paper Oka's conjecture on irreducible plane sextics
\jour J. London Math. Soc.
\toappear
\finalinfo\tt arXiv:\allowbreak math.AG/0701671
\endref\label{degt.Oka}

\ref{D5}
\by A.~Degtyarev
\paper Zariski $k$-plets via dessins d'enfants
\jour Comment. Math. Helv.
\finalinfo\tt arXiv:0710.0279
\toappear
\endref\label{degt.kplets}

\ref{D5}
\by A.~Degtyarev
\paper On irreducible sextics with non-abelian fundamental group
\inbook Fourth Franco-Japanese Symposium on Singularities (Toyama, 2007)
\finalinfo\tt arXiv:0711.3070
\toappear
\endref\label{degt.Oka3}

\ref{D6}
\by A.~Degtyarev
\paper Irreducible plane sextics with large fundamental groups
\finalinfo\tt arXiv:0712.2290
\endref\label{degt.8a2}

\ref{D7}
\by A.~Degtyarev
\paper Stable symmetries of plane sextics
\finalinfo\tt arXiv:0802.2336
\endref\label{symmetric}

\ref{EO}
\by C.~Eyral, M.~Oka
\paper Fundamental groups of the complements of certain plane
non-tame torus sextics
\jour Topology Appl.
\vol 153
\yr 2006
\issue 11
\pages 1705--1721
\endref\label{EyralOka}

\ref{GAP}
\by The GAP Group
\book GAP --- Groups, Algorithms, and Programming
\bookinfo Version 4.4.10
\yr 2007
\finalinfo ({\tt http://www.gap-system.org})
\endref\label{GAP}

\ref{vK}
\by E.~R.~van~Kampen
\paper On the fundamental group of an algebraic curve
\jour  Amer. J. Math.
\vol   55
\yr    1933
\pages 255--260
\endref\label{vanKampen}

\ref{No}
\by M.~V.~Nori
\paper Zariski conjecture and related problems
\jour Ann. Sci. \'Ec. Norm. Sup., 4 s\'erie
\vol    16
\yr     1983
\pages  305--344
\endref\label{Nori}

\ref{OP1}
\by M.~Oka, D.~T.~Pho
\paper Classification of sextics of torus type
\jour Tokyo J. Math.
\vol 25
\issue 2
\pages 399--433
\yr 2002
\endref\label{OkaPho.moduli}

\ref{OP2}
\by M.~Oka, D.~T.~Pho
\paper Fundamental group of sextics of torus type
\inbook Trends in singularities
\pages 151--180
\bookinfo Trends Math.
\publ Birkh\"auser
\publaddr Basel
\yr 2002
\endref\label{OkaPho}

\ref{Oz}
\by A.~\"Ozg\"uner
\book Classical Zariski pairs with nodes
\bookinfo M.Sc. thesis
\publ Bilkent University
\yr 2007
\endref\label{Aysegul}

\ref{Tju}
\by G.~Tjurina
\paper Resolution of singularities of flat deformations of double rational
points
\jour Funktsional Anal. i Pril.
\vol 4
\yr 1970
\issue 1
\pages 77--83
\lang Russian
\transl\nofrills English transl. in
\jour Functional Anal. Appl.
\vol 4
\issue 1
\yr 1970
\pages 68-73
\endref\label{Tjurina}

\ref{T}
\by H.~Tokunaga
\paper $(2,3)$-torus sextics and the Albanese images of $6$-fold
cyclic multiple planes
\jour Kodai Math.~J.
\vol 22
\yr 1999
\issue 2
\pages 222--242
\endref\label{Tokunaga}

\ref{Z}
\by O.~Zariski
\paper On the problem of existence of algebraic functions of two
variables possessing a given branch curve
\jour Amer. J. Math.
\vol 51
\yr 1929
\pages 305--328
\endref\label{Zariski}

\ref{Ya}
\by J.-G.~Yang
\paper Sextic curves with simple singularities
\jour Tohoku Math. J. (2)
\vol 48
\issue 2
\yr 1996
\pages 203--227
\endref\label{Yang}

\endRefs
\enddocument